\documentclass[12pt]{amsart}
\usepackage{amsmath}
\usepackage{amsxtra}
\usepackage{amscd}
\usepackage{amsthm}
\usepackage{amsfonts}
\usepackage{amssymb}
\usepackage{eucal}
\usepackage{mathrsfs}
\usepackage{color, graphics}
\usepackage{hyperref}

\def\C{{\mathbb C}}

\def\N{{\mathbb N}}
\def\P{{\mathbb P}}

\def\Z{{\mathbb Z}}

\newtheorem{theorem}{Theorem}[section]
\newtheorem{definition}{Definition}[section]
\newtheorem{lemma}[theorem]{Lemma}

\newtheorem{proposition}[theorem]{Proposition}
\newtheorem{corollary}[theorem]{Corollary}

\title{On Higher Syzygies of Ruled varieties over a Curve}

\author[E. Park]{Euisung Park}
\thanks{Supported by Korea Research Foundation Grant (KRF-2002-070-C00003)}
\address {Euisung Park : School of Mathematics, Korea Institute for Advanced Study,
270-43 Cheongryangri-dong, Dongdaemun-gu, Seoul 130-012, Republic of Korea,} \email{puserdos@kias.re.kr}

\begin{document}

\thispagestyle{empty} \maketitle

\begin{abstract}
For a vector bundle $\mathcal{E}$ of rank $n+1$ over a smooth
projective curve $C$ of genus $g$, let $X=\P_C (\mathcal{E})$ with
projection map $\pi:X\rightarrow C$. In this paper we investigate
the minimal free resolution of homogeneous coordinate rings of
$X$.

We first clarify the relations between higher syzygies of very
ample line bundles on $X$ and higher syzygies of Veronese
embedding of fibres of $\pi$ by the same line bundle. More
precisely, letting $H = \mathcal{O}_{\P_C (\mathcal{E})} (1)$ be
the tautological line bundle, we prove that if
$(\P^n,\mathcal{O}_{\P^n} (a))$ satisfies Property $N_p$, then
$(X,aH+\pi^*B)$ satisfies Property $N_p$ for all $B \in
\mbox{Pic}C$ having sufficiently large degree(Theorem
\ref{thm:positive}). And also the effective bound of
$\mbox{deg}(B)$ for Property $N_p$ is obtained(Theorem
\ref{thm:1}, \ref{thm:2}, \ref{thm:3} and \ref{thm:4}). For the
converse, we get some partial answer(Corollary
\ref{cor:negative}). Secondly, by using these results we prove
some Mukai-type statements. In particular, Mukai's conjecture is
true for $X$ when $\mbox{rank}(\mathcal{E}) \geq g$ and $\mu^-
(\mathcal{E})$ is an integer(Corollary \ref{cor:Mukai}). Finally
for all $n$, we construct an $n$-dimensional ruled variety $X$ and
an ample line bundle $A \in \mbox{Pic}X$ which shows that the
condition of Mukai's conjecture is optimal for every $p \geq 0$.
\end{abstract}

\tableofcontents \setcounter{page}{1}

\section{Introduction} \noindent For a projective variety $X$ and
a very ample line bundle $L$ on $X$, the equations defining $X
\subset \P = \P H^0 (X,L)$ and the higher syzygies among them have
been studied by several authors(\cite{Butler}, \cite{EL},
\cite{Green}, \cite{Green2}, \cite{GL}, \cite{OP}, \cite{Ru},
etc). Classical results is about the numerical and cohomological
conditions to guarantee that $X$ is projectively normal or cut out
by quadrics. Nowadays this problem is generalized to higher
syzygies. We recall the notion of Property $N_p$ introduced by M.
Green and R. Lazarsfeld\cite{GL}. Let $S$ be the homogeneous
coordinate ring of $\P$, $I_X \subseteq S$ the homogeneous ideal
of $X$, and $S(X) = S/I_X$ the homogeneous coordinate ring of $X$.
For the graded $S$-module $E:= \oplus_{\ell \in \Z} H^0
(X,L^{\ell})$, let us consider the minimal free resolution:
\begin{equation*}
\cdots \rightarrow L_i
\rightarrow \cdots \rightarrow L_1
\rightarrow L_0 \rightarrow E \rightarrow 0
\end{equation*}
where $L_i$, as free graded $S$-module, can be written as
\begin{equation*}
L_i = \oplus_j S^{k_{i,j}}(-i-j).
\end{equation*}
M. Green and R. Lazarsfeld\cite{GL}
considered the situations in which the first few modules of
syzygies of $E$ are as simple as possible:\\

\begin{definition} For a nonnegative integer $p$, $L$ is said
to satisfy Property $N_p$ if $k_{i,j} = 0$ for $0 \leq i \leq p$
and  $j \geq 2$. Equivalently, Property $N_p$ holds for $X
\hookrightarrow \P$ if $E$ admits a minimal free resolution of the
form
\begin{eqnarray*}
\cdots \rightarrow S^{k_{p,1}}(-p-1) \rightarrow \cdots
\rightarrow S^{k_{2,1}}(-3) \rightarrow  S^{k_{1,1}}(-2)
\rightarrow S \rightarrow E \rightarrow 0.\\
\end{eqnarray*}
\end{definition}

\noindent So, remark that $L$ satisfies Property $N_0$ if and only
if $L$ is normally generated(i.e., projectively normal), $L$
satisfies Property $N_1$ if and only if $L$ satisfies Property
$N_0$ and the homogeneous ideal is generated by quadrics, $L$
satisfies Property $N_2$ if and only if $L$ satisfies Property
$N_1$ and the relations among the quadrics are generated
by the linear relations and so on.\\

\noindent {\bf Remark 1.1.} Given a projective variety $X$ and a
very ample line bundle $L \in \mbox{Pic}X$, it was expected that
if $L$ is more and more ``positive" then the embedding $X
\hookrightarrow \P H^0 (X,L)$ is more and more ``regular". And
several results indicate that the ``regularity" should be Property
$N_p$. See Theorem 3.2 in \cite{Green2}. \qed \\

\noindent {\bf Remark 1.2.} Let $L$ be a very ample line bundle
with $H^1 (X,L^j)=0$ for all $j \geq 2$. Recently Sijong Kwak and
the author proved in \cite{KP} that if $L$ satisfies Property
$N_p$, then every embedding given by a very ample subsystem in
$H^0 (X,L)$ of codimension $c \leq p-1$ satisfies $k$-normality
for all $k \geq c+1$ and the homogeneous ideal is generated by
forms of degree $\leq c+2$. In this sense, Property $N_p$ is also
very important in classical projective geometry. \qed \\

\noindent The first effective result in this field is M. Green's
theorem for curves. In \cite{Green}, he proved that when $C$ is a
smooth curve of genus $g$ and $\deg L \geq 2g+1+p~$, then $(C,L)$
satisfies Property $N_p$. Then Shigeru Mukai suggests what form a
statement should take for higher dimensional varieties. He
observed that one could view Green's Theorem as asserting that if
$D$ is an ample bundle on $C$, then Property $N_p$ holds for $K_C
+(p+3)D$. This leads him to conjecture that for a smooth
projective variety $X$ of dimension $n$ and an ample line bundle
$A \in \mbox{Pic}X$, $K_X + (n+2+p)A$ satisfies Property $N_p$.
Nowadays his conjecture is a popular guiding principle for
studying higher syzygies of arbitrary varieties and many related
remarkable research have been
developed. \\

In this paper we study Property $N_p$ of very ample line bundles
on ruled varieties over a curve. Although Mukai's conjecture
suggests a guideline to study higher syzygies for arbitrary
varieties, it should be a more natural question that for a given
very ample line bundle $L$ on a smooth projective variety $X$,
determine the largest $p$ such that Property $N_p$ holds for $L$.
And our first result gives a partial answer for this question when
$X$ is a ruled variety over a curve(Theorem \ref{thm:positive}).
Let $C$ be a smooth projective curve of genus $g$. For a vector
bundle $\mathcal{E}$ of rank $n+1$ over $C$, let $X=\P_C
(\mathcal{E})$ be the associated projective space bundle with
tautological line bundle $H=\mathcal{O}_{\P_C (\mathcal{E})} (1)$
and projection map $\pi:X\rightarrow C$. Therefore
\begin{equation*}
\mbox{Pic}X = \Z[H]\oplus \pi^* \mbox{Pic}C
\end{equation*}
and every line bundle is written by $aH+\pi^* B$ for some $a \in
\Z$ and $B \in \mbox{Pic}C$. If $L = aH+\pi^*B$ is a very ample
line bundle, then the embedding $X \hookrightarrow \P H^0 (X,L)$
is a fiberwise $a$-uple Veronese embedding. Thus it seems natural
to investigate algebraic relations between the $a$-uple Veronese
embedding $\P^n \hookrightarrow \P H^0 (\P^n,\mathcal{O}_{\P^n}
(a))$ and the embedding $X \hookrightarrow \P H^0 (X,L)$. Along
this point of view we prove

\begin{theorem}\label{thm:positive}
Assume that $(\P^n,\mathcal{O}_{\P^n} (a))$ satisfies Property
$N_p$ for some $p \geq 0$. Then there is an integer $s(a,p)$ such
that $(X,aH+\pi^*B)$ satisfies Property $N_p$ for any $B \in
\mbox{Pic}C$ such that $\mbox{deg}(B) \geq s(a,p)$.
\end{theorem}

\noindent The key idea is that if the fibre
$(\P^n,\mathcal{O}_{\P^n} (a))$ satisfies Property $N_p$, then
Property $N_p$ of $aH+\pi^* B$ can be checked by the vanishing of
some cohomology groups on $C$. For the converse it may be expected
that if $(\P^n,\mathcal{O}_{\P^n}(a))$ fails to satisfy Property
$N_p$ then $(X,aH+\pi^*B)$ also fails to satisfy Property $N_p$ no
matter how large $\mbox{deg}(B)$ is. At least for $C \times \P^n$
(i.e. $\mathcal{E} = \oplus^{n+1} \mathcal{O}_C$), this is true by
Elena Rubei's result\cite{Ru}. For arbitrary $\mathcal{E}$ we have
a partial answer. First we recall known results about higher
syzygies of Veronese embedding. For details, see \cite{OP}.\\

\begin{enumerate}
\item[$(1)$] $(\P^n,\mathcal{O}_{\P^n}(d))$  satisfies Property
$N_d$. \item[$(2)$] $(\P^n,\mathcal{O}_{\P^n}(d))$ satisfies
Property $N_p$ for all $p \geq 0$ if and only if $n=1$ or $a=1$ or
$n=a=2$. \item[$(3)$] For $n \geq 3$,
$(\P^n,\mathcal{O}_{\P^n}(2))$  satisfies Property $N_p$ if and
only if $p \leq 5$. \item[$(4)$] For $d \geq 3$,
$(\P^2,\mathcal{O}_{\P^n}(d))$  satisfies Property $N_p$ if and
only if $p \leq 3d-3$. \item[$(5)$] $(\P^3,\mathcal{O}_{\P^3}(3))$
satisfies Property $N_{p}$ if and only if $p \leq 6$. \item[$(6)$]
For $n \geq 2$ and $d \geq 3$, $(\P^n,\mathcal{O}_{\P^n}(d))$
fails to satisfy Property $N_{3d-2}$.
\end{enumerate}

\begin{table}[hbt]
\begin{center}
\begin{tabular}{|c|c|c|c|c| }\hline
$a \backslash n$ & \quad $1$ \quad                 & $\quad \quad
2 \quad \quad $  & $3  \cdots \cdots$      \\\hline
 $ 1 $           &  \multicolumn{3}{c|}{}                                                  \\\cline{1-1} \cline{4-4}
 $ 2 $           & \multicolumn{2}{c|}{Property $N_p$ holds for all $p \geq 0$.} & Property $N_p$ holds if and only if $p\leq5$.   \\\cline{1-1} \cline{3-4}
 $ 3 $           &   &              & Property $N_a$ holds. \\\cline{1-1}
 $ 4 $           &   &  Property $N_p$ holds          &  Property $N_p$ fails to hold         \\
$\vdots$         &                      &  if and only if $p \leq
3d-3$. &  for all $p \geq 3d-2$.  \\\hline
\end{tabular}
\end{center}
\caption{Property $N_p$ for $(\P^n,\mathcal{O}_{\P^n}(d))$}
\end{table}

\noindent See \cite{Ru2} for the proof that
$(\P^n,\mathcal{O}_{\P^n}(3))$ satisfies Property $N_{4}$. Also G.
Ottaviani and R. Paoletti conjectured that for $n \geq 2$ and $d
\geq 3$, $(\P^n,\mathcal{O}_{\P^n}(d))$ satisfies Property $N_p$
if and only if $p \leq 3d-3$\cite{OP}. As an analogue of $(3)$ and
$(6)$, we prove a partial converse of Theorem \ref{thm:positive}.
Indeed it is induced immediately from D. Eisenbud, M. Green, K.
Hulek and S. Popescu's recent work in \cite{EGHP}. See Theorem
\ref{thm:EGHP} and Corollary \ref{cor:negative}.

By Theorem \ref{thm:positive}, if $(\P^n,\mathcal{O}_{\P^n}(a))$
satisfies Property $N_p$ then it is natural to find a bound of
$s(a,p)$. And Theorem \ref{thm:1}(ruled scrolls), Theorem
\ref{thm:2}(ruled surfaces), Theorem \ref{thm:3}(Veronese surface
fibrations) and Theorem \ref{thm:4}(arbitrary cases) are devoted
to compute a concrete bound of $s(a,p)$. Also see Corollary
\ref{cor:detailscroll}, Corollary \ref{cor:detailsurface} and
Corollary \ref{cor:detailarbitrary}. Our bounds are optimal when
$g \leq 2$ as discussed in Example 4.1.1 and Example 4.4.1. The
results are stated in terms of the minimal slope $\mu ^-(p_* L)$
of the vector bundle $\pi_* L$. The minimal slope is a natural
invariant to measure the positivity of line bundles on $X$. Thus
Theorem \ref{thm:1}, Theorem \ref{thm:2} and Theorem \ref{thm:3}
imply that when $a=1$ or $n=1$ or $n=a=2$(equivalently when
$(\P^n,\mathcal{O}_{\P^n} (a))$ satisfies Property $N_p$ for all
$p \geq 0$), the minimal free resolution of the embedding $X
\hookrightarrow \P H^0 (X,L)$ is more and more simple as the
vector bundle $\pi_* L$ is more and more ``positive". Also by
applying Theorem \ref{thm:4} we  obtain a Mukai type result for
$X$. See Corollary \ref{cor:Mukai}. In particular we prove that if
$\mbox{rank}(\mathcal{E}) \geq g$ and $\mu^- (\mathcal{E})$ is an
integer, then Mukai's conjecture is true for $X$.

Although Mukai's conjecture is still open even for ruled varieties
over a curve, it seems worthwhile to consider the optimality of
the statement. With this end in view, we construct some
interesting examples. Indeed for all $n \geq 2$ and $p \geq 0$ we
construct an $n$-dimensional ruled variety $X$ over a curve and an
ample line bundle $A \in \mbox{Pic}X$ such that for any $p \geq
0$, $(X,K_X + (n+2+p)A)$ fails to satisfy Property $N_{p+1}$.

We should point out that higher syzygies of ruled varieties over a
curve are already investigated by David C. Butler\cite{Butler}.

\begin{theorem}[David C. Butler,
\cite{Butler}]\label{thm:Butlermain} For a rank $n$ vector bundle
$\mathcal{E}$  over $C$ a smooth projective curve of genus $g$,
let $X=\P(\mathcal{E})$ with tautological line bundle $H$ and
projection map $\pi:X \rightarrow C$. Put $L=aH+\pi^* B$ where
$a\geq1$ and $B\in PicC)$.
\begin{enumerate}
\item[$(i)$] If $\mu^- (\pi_* L)\geq 2g+1$, then
$(X,L)$ is normally generated.
\item[$(ii)$] Fix an integer $1\leq p\leq a-1$. If $\mu^- (\pi_* L)\geq 2g+2p$,
then Property $N_p$ holds for $(X,L)$.
\item[$(iii)$] If $A_i$ is an ample line bundles on $X$, then $(X,K_X
+ A_1 + \cdots + A_t)$ satisfies Property $N_p$ for all $t \geq
2n+2np$.
\end{enumerate}
\end{theorem}

The organization of this paper is as follows. In $\S 2$, we review
some necessary elementary facts about vector bundles over a curve,
regularity over ruled varieties over a variety of arbitrary
dimension, Bott formula, etc. In $\S 3$, we explain that the
relation between higher syzygies of very ample line bundles on
ruled varieties and those of the Veronese embedding as fibres of
$\pi$. Then we prove Theorem \ref{thm:positive} and Corollary
\ref{cor:negative}. $\S 4$ is devoted to prove some results about
Property $N_p$ of ruled varieties over a curve. Finally in $\S 5$
we produce some interesting examples which shows that the
numerical condition of Mukai's conjecture is optimal.\\

\section{Notation and Preliminaries}
\subsection{Notations} Throughout this paper the following is assumed.
\begin{enumerate}
\item[$1.$] All varieties are defined over the complex number
field $\C$.
\item[$2.$] For a finite dimensional $\C$-vector space $V$, $\P(V)$ is
the projective space of one-dimensional quotients of $V$.
\item[$3.$] When a projective variety $X$ is embedded in a projective space
$\P^r$ by a very ample line bundle $L$ on it, we may write
$\mathcal{O}_X (1)$ instead of $L$ so long as no confusion arise.
\end{enumerate}

\subsection{The slope of vector bundles over a curve }
Let $C$ be a smooth projective curve of genus $g$. For a vector
bundle $\mathcal{F}$ on $C$, the slope $\mu(\mathcal{F})$ is
defined by $\deg(\mathcal{F}) / \mbox{rank}(\mathcal{F})$. Also
the maximal slope and the minimal slope are defined as follows:
\begin{equation*}
\mu^+ (\mathcal{F})= \mbox{max}\{ \mu(S) | 0 \rightarrow S
\rightarrow \mathcal{F}\} ~~~\mbox{and}~~~ \mu^-
(\mathcal{F})=\mbox{min}\{ \mu(Q) | \mathcal{F}\rightarrow Q
\rightarrow 0\}
\end{equation*}
These notions  satisfy the following properties.

\begin{lemma}\label{lem:folklore}
For vector bundles $\mathcal{E}$, $\mathcal{F}$ and $\mathcal{G}$ on $C$,\\
$(1)$ $\mu^+ ( \mathcal{E} \otimes \mathcal{F}) = \mu^+
(\mathcal{E} ) +\mu^+ ( \mathcal{F})$.\\
$(2)$ $\mu^- ( \mathcal{E} \otimes \mathcal{F}) = \mu^-
(\mathcal{E} ) +\mu^- ( \mathcal{F})$.\\
$(3)$ $\mu^+ (S^\ell (\mathcal{E})) = \ell \mu^+ (\mathcal{E})$.\\
$(4)$ $\mu^- (S^\ell (\mathcal{E})) = \ell \mu^- (\mathcal{E})$.\\
$(5)$ $\mu^- (\wedge^\ell \mathcal{E}) \geq \ell \mu^-
(\mathcal{E})$.\\
$(6)$ If $\mu^- (\mathcal{E}) > 2g-2$, then $h^1 (C, \mathcal{E})=0$.\\
$(7)$ If $\mu^- (\mathcal{E})  > 2g-1$, then $\mathcal{E}$ is globally generated.\\
$(8)$ If $\mu^- (\mathcal{E})  > 2g$, then
$\mathcal{O}_{\P(\mathcal{E})} (1)$ is very ample.\\
$(9)$ If $0 \rightarrow \mathcal{E} \rightarrow \mathcal{F}
\rightarrow \mathcal{G} \rightarrow 0$ is an exact sequence, then
\begin{eqnarray*}
\mu^- (\mathcal{G}) \geq \mu^- (\mathcal{F}) \geq \mbox{min}\{
\mu^- (\mathcal{E}), \mu^- (\mathcal{G}) \}.\\
\end{eqnarray*}
\end{lemma}

\begin{proof}
See $\S 1$ and $\S 2$ in \cite{Butler}.
\end{proof}

\noindent Therefore  if $\mu^- (\mathcal{E}) > 2g-1$, then the
evaluation map determines the short exact sequence
\begin{equation*}
0 \rightarrow M_{\mathcal{E}} \rightarrow H^0 (C,\mathcal{E})
\otimes \mathcal{O}_C \rightarrow \mathcal{E} \rightarrow 0.
\end{equation*}

\noindent And David C. Butler obtained the following very useful
result:

\begin{theorem}[David C. Butler, \cite{Butler}]\label{thm:Butlerestimation}
For a vector bundle $\mathcal{E}$ over $C$, if $\mu^-
(\mathcal{E})\geq 2g$, then $M_{\mathcal{E}}$ satisfies
\begin{equation*}
\mu^- (M_{\mathcal{E}}) \geq -\frac{\mu^- (\mathcal{E})}{\mu^-
(\mathcal{E})-g}.
\end{equation*}\\
\end{theorem}

\noindent {\bf Remark 3.1.} Clearly $\mu^- (M_{\mathcal{E}}) \leq
\mu (M_{\mathcal{E}}) < 0$. Indeed, since $\det (M_{\mathcal{E}})
= - \det(\mathcal{E}) <0$,
\begin{equation*}
\mu (M_{\mathcal{E}}) = \det (M_{\mathcal{E}}) / \mbox{rank}
(M_{\mathcal{E}}) <0.
\end{equation*}
This fact will be used in several inequalities. \qed \\

Finally we recall Miyaoka's criterion for ampleness of line
bundles on ruled varieties over a curve.

\begin{lemma}[Y. Miyaoka]\label{lem:ampleMiyaoka}
Let $\mathcal{E}$ be a vector bundle on $C$, and let $X = \P_C
(\mathcal{E})$ with $H= \mathcal{O}_{\P_C (\mathcal{E})} (1)$.
Then $aH+\pi^* B \in \mbox{Pic}X$ is ample if and only if $a \geq
1$ and $a \mu^- (\mathcal{E}) + \mbox{deg} (B) >0$.
\end{lemma}

\begin{proof}
See Lemma 5.3 in \cite{Butler}.
\end{proof}

\subsection{Regularity of vector bundles over ruled varieties}
We recall some basic facts about the regularity of vector bundles
over ruled varieties. Let $E$ be a vector bundle of rank $n+1$
over a smooth projective variety $Y$ and let $X=\P_Y (E)$ with the
projection map $\pi : X \rightarrow Y$ and tautological line
bundle $H$.

\begin{definition}
For a vector bundle $\mathcal{F}$ over $X$, we say that
$\mathcal{F}$ is $f$  $\pi$-regular when
\begin{equation*}
R^i \pi_* (\mathcal{F}(f-i)) = 0
\end{equation*}
for every $i \geq 1$.
\end{definition}

\noindent Here $\mathcal{F}(f-i) =\mathcal{F} \otimes (f-i)H$.
Note that a line bundle of the form $aH+\pi^* B$ is $(-a)$
$\pi$-regular. We present some basic facts about the
$\pi$-regularity.

\begin{lemma}\label{lem:directimage}
Let $\mathcal{F}$ and $\mathcal{G}$ be two vector bundles on $X$
with $f$ and $g$ $\pi$-regularity, respectively.\\
$(1)$ $\mathcal{F} \otimes \mathcal{G}$ is $(f+g)$
$\pi$-regular.\\
$(2)$ If $f\leq 1$, then
\begin{equation*}
H^i (X,\mathcal{F})\cong H^i (Y,\pi_* \mathcal{F})\quad \mbox{for
all}\quad i\geq0.
\end{equation*}
$(3)$ If $f\leq0$ and $\widetilde{\mathcal{F}}=\pi^* (\pi_*
\mathcal{F})$, there is an exact sequence of vector bundles on $X$
\begin{equation*}
0 \rightarrow  \mathcal{K}_{\widetilde{\mathcal{F}}} \rightarrow
\widetilde{\mathcal{F}}  \rightarrow \mathcal{F} \rightarrow 0
\end{equation*}

where $\mathcal{K}_{\widetilde{\mathcal{F}}}$ is $1$ $\pi$-regular.\\
$(4)$ If $f \leq 0$ and $g \leq 0$, then there is a surjective map
\begin{equation*}
\pi_* \mathcal{F} \otimes \pi_*  \mathcal{G} \rightarrow \pi_*
(\mathcal{F} \otimes \mathcal{G}) \rightarrow 0.
\end{equation*}

In particular, if $Y$ is a curve
\begin{equation*}
\mu^- (\pi_* (\mathcal{F} \otimes \mathcal{G})) \geq \mu^-  (\pi_*
\mathcal{F}) +  \mu^- (\pi_*  \mathcal{G}).
\end{equation*}
\end{lemma}

\begin{proof}
See Lemma 3.2 in \cite{Butler}.
\end{proof}

\subsection{Bott formula}
Let $Y$ be a projective variety and let $\mathcal{E}$ be a vector
bundle of rank $n+1$. Let $X=\P_Y (\mathcal{E})$ with the
tautological line bundle $\mathcal{O}_X (1)$ and projection
morphism $\pi : X \rightarrow Y$. Then there is a natural exact
sequence
\begin{equation*}
0 \rightarrow \Omega_{X/Y} (1) \rightarrow \pi^* \mathcal{E}
\rightarrow \mathcal{O}_X (1) \rightarrow 0
\end{equation*}
where $\Omega_{X/Y}$ denotes the relative canonical sheaf which is
clearly of rank $n$. Denote $\wedge^j \Omega_{X/Y}$ by
$~\Omega_{X/Y} ^j$.

\begin{proposition}\label{prop:Bott}
$(1)$ For $1\leq j \leq n$ and $k >j$,
\begin{equation*}
\pi_*  \Omega_{X/Y} ^j (k) ~~~~\mbox{ is a vector bundle on $Y$ of
rank }~~~~ {{k+n-j}\choose{k}} {{k-1}\choose{j}}.
\end{equation*}
$(2)$ For $i \geq 1,~1\leq j \leq n$ and $k \geq 1$,
\begin{equation*}
R^i \pi_*  \Omega_{X/Y} ^j (k) = 0.
\end{equation*}
$(3)$ Let $\mathcal{F}$ be a vector bundle on $Y$. Then
\begin{equation*}
H^i (X,\Omega_{X/Y} ^j (k) \otimes \pi^* \mathcal{F}) \cong H^i
(Y,\pi_*  \Omega_{X/Y} ^j (k) \otimes \mathcal{F})
\end{equation*}
for $i \geq 1,~1\leq j \leq n$ and $k \geq 1$.
\end{proposition}

\begin{proof}
For every $y \in Y$,
\begin{equation*}
\begin{CD}
R^i \pi_*  \Omega_{X/Y} ^j (k) |_{\pi^{-1} (y)} & \cong & H^i (\pi^{-1} (y), \Omega_{X/Y} ^j (k) |_{\pi^{-1} (y)}) \\
                                             & \cong & H^i (\P^{n},\Omega_{\P^{n}} ^j (k) ).
\end{CD}
\end{equation*}
Then by Bott formula,
$$h^i (\P^{n},\Omega_{\P^{n}} ^j (k) ) = \begin{cases} {{k+n-j}\choose{k}}
{{k-1}\choose{j}}
~~~~& \mbox{if $i=0,1\leq j \leq n$ and $k >j$, and}\\
0 ~~~~& \mbox{if $i \geq 1,~1\leq j \leq n$ and $k \geq 1$}
\end{cases}$$\\
which completes the proof of $(1)$ and $(2)$. Finally $(3)$
follows from
\begin{equation*}
R^i \pi_* (\Omega_{X/Y} ^j (k) \otimes \pi^* \mathcal{F}) = R^i
\pi_* (\Omega_{X/Y} ^j (k)) \otimes \mathcal{F}=0
\end{equation*}
for $i \geq 1,~1\leq j \leq n$ and $k \geq 1$ which is guaranteed
by $(2)$ and projection formula.
\end{proof}

\subsection{Criteria for Property $N_p$} In this subsection we
recall well-known cohomological criteria for Property $N_p$. Let
$X$ be a smooth projective variety and let $L \in \mbox{Pic}X$ be
a very ample line bundle such that $H^1 (X,L^j )=0$ for all
$j\geq1$. Consider the short exact sequence
\begin{equation*}
0 \rightarrow M_L \rightarrow  H^0 (X,L) \otimes \mathcal{O}_X
\rightarrow L \rightarrow 0
\end{equation*}
where $H^0 (X,L) \otimes \mathcal{O}_X \rightarrow L$ is the
natural surjective evaluation map.

\begin{lemma}\label{lem:criterion}
$(1)$ $(X,L)$ satisfies Property $N_p$ if and only if $H^1
(X,\wedge^i M_L \otimes L^j )=0$

for all $1\leq i \leq p+1$ and $j\geq1$.\\
$(2)$ For $p \leq \mbox{codim}(X,\P H^0 (X,L))$, Property $N_p$
holds for $(X,L)$ if and only if
\begin{equation*}
H^1 (X,\wedge^{p+1} M_L \otimes L^j )=0~~\mbox{for all $j\geq1$.}
\end{equation*}
\end{lemma}

\begin{proof}
See $\S 1$ in \cite{GL2} and Lemma 1.6 in \cite{EL}.
\end{proof}

\section{Veronese Embedding and Property $N_p$ of
Ruled Varieties} \noindent Let $C$ be a smooth projective curve of
genus $g$ and let $\mathcal{E}$ be a vector bundle of rank $n+1$
over $C$. For the associated projective space bundle $X=\P_C
(\mathcal{E})$, let $H = \mathcal{O}_{\P_C (\mathcal{E})} (1)$ be
the tautological line bundle and $\pi : X \rightarrow C$ the
projection morphism. Thus
\begin{equation*}
\mbox{Pic}X = \Z[H]\oplus \pi^* \mbox{Pic}C
\end{equation*}
and  every line bundle on $X$ is written as $aH+\pi^* B$ for some
$a\in \Z$ and $B\in \mbox{Pic}C$.\\

\begin{lemma}\label{lem:veryampleness}
For $L=aH+\pi^*B \in \mbox{Pic}X$, assume that $\mu^- (\pi_* L)  > 2g$. Then\\
$(1)$ $L$ is very ample.\\
$(2)$ $H^1 (X,L^j)=0$ for all $j \geq 1$.\\
$(3)$ $\mbox{codim} (X,\P H^0 (X,L)) \geq \begin{cases} 5
~~~~& \mbox{if  $n\geq3$ and $a=2$, and}\\
3a-3~~~~& \mbox{if $n \geq 2$ and $a \geq 3$}.
\end{cases}$
\end{lemma}

\begin{proof}
$(1)$ When $a=1$, see Lemma \ref{lem:folklore}.$(8)$. When
$a\geq2$, consider the vector bundle $\mathcal{F}=\pi_* L$. Then
$X \subset \P_C (\mathcal{F})$ is given by a fiberwise $a$-uple
map and
\begin{equation*}
\mathcal{O}_{\P_C (\mathcal{F})} (1)|_{X} = L.
\end{equation*}
Since $\mathcal{O}_{\P_C (\mathcal{F})} (1)$ is very ample, so is
$L$.\\
$(2)$ See Lemma \ref{lem:folklore}.$(6)$.\\
$(3)$ Note that $h^0 (X,L) = h^0 (C,\mathcal{F})$ where
$\mathcal{F}=\pi_* L$. Also note that
\begin{equation*}
\mbox{rank}(\mathcal{F})={{n+a}\choose{a}}~~\mbox{and}~~
\mbox{degree}(\mathcal{F})={{n+a}\choose{a-1}} \cdot \mbox{deg}
\mathcal{F} + \mbox{rank}(\mathcal{F}) \cdot \mbox{deg}(B).
\end{equation*}
Thus
\begin{equation*}
h^0 (C,\mathcal{F})\geq {{n+a}\choose{a}} (\mu^- (\mathcal{F}) -
g+1)  \geq {{n+a}\choose{a}}
\end{equation*}
by Riemann-Roch. Therefore it suffices to show that
\begin{equation*}
{{n+a}\choose{a}}-(n+1) \geq \begin{cases} 5
~~~~& \mbox{if  $n \geq 3$ and $a=2$, and}\\
3a-3~~~~& \mbox{if $n \geq 2$ and $a \geq 3$}
\end{cases}
\end{equation*}
which is easily checked.
\end{proof}

\noindent {\bf Proof of Theorem \ref{thm:positive}.} Since $\mu^-
(\pi_* L)  = a \mu^- (\mathcal{E}) + \mbox{deg}(B)$, we may assume
that $\mu^- (\pi_* L) > 2g$. Thus $L$ is very ample and $H^1
(X,L^j)=0$ for all $j \geq 1$ by Lemma \ref{lem:veryampleness}.
Consider the short exact sequence
\begin{equation*}
0 \rightarrow M_L \rightarrow  H^0 (X,L) \otimes \mathcal{O}_X
\rightarrow L \rightarrow 0.
\end{equation*}
We first prove that $M_L$ is $1$ $\pi$-regular. Let
$\mathcal{G}_a$ be the vector bundle on $X$ defined by the short
exact sequence
\begin{equation*}
0 \rightarrow \mathcal{G}_a \rightarrow \pi^* S^a (\mathcal{E})
\rightarrow aH \rightarrow 0.
\end{equation*}
Note that for each $P \in C$, the restriction of the above
sequence to $\pi^{-1} (P) \cong \P^n$ is just equal to the short
exact sequence
\begin{equation*}
0 \rightarrow M_a \rightarrow H^0 (\P^n,\mathcal{O}_{\P^n}
(a))\otimes \mathcal{O}_{\P^n}  \rightarrow  \mathcal{O}_{\P^n}
(a) \rightarrow 0.
\end{equation*}
Since we assume that $(\P^n,\mathcal{O}_{\P^n} (a))$ satisfies
Property $N_p$, $H^i (\P^n,\wedge^m M_a \otimes \mathcal{O}_{\P^n}
(aj))=0$ for $1 \leq m \leq p+1$ and all $i,j \geq 1$. Thus
\begin{equation*}
R^i \pi_* (\wedge^m \mathcal{G}_a \otimes ajH) =0~~\mbox{for $1
\leq m \leq p+1$ and all $i,j \geq 1$}. \quad \quad (*)
\end{equation*}
Also $\mathcal{G}_a$ is $1$ $\pi$-regular by Lemma
\ref{lem:directimage}.$(3)$. Put $\mathcal{F} = \pi_* L =S^a
(\mathcal{E}) \otimes B$. Then we have the following commutative
diagram:
\begin{equation*}
\begin{CD}
&  &&  && 0 &\\
&  &&  && \downarrow &\\
& 0 && && \mathcal{G}_a \otimes \pi^* B &\\
& \downarrow &&  && \downarrow &\\
0 \rightarrow &\pi^* M_\mathcal{F}& \rightarrow & H^0 (C,\mathcal{F}) \otimes \mathcal{O}_X & \rightarrow  & \pi^* \mathcal{F} & \rightarrow 0 \\
& \downarrow && \parallel && \downarrow &\\
0 \rightarrow & M_L & \rightarrow & H^0 (X,L) \otimes \mathcal{O}_X & \rightarrow & L & \rightarrow 0 \\
& \downarrow &&  && \downarrow &\\
& \mathcal{G}_a \otimes \pi^* B && && 0  &\\
& \downarrow  &&  && &\\
& 0 &&  && &.\\
\end{CD}
\end{equation*}
Since $\pi^* M_\mathcal{F}$ is $0$ $\pi$-regular, it is
immediately proved that $M_L$ is $1$ $\pi$-regular from the short
exact sequence $0 \rightarrow \pi^* M_\mathcal{F} \rightarrow M_L
\rightarrow \mathcal{G}_a \otimes \pi^* B \rightarrow 0$.
Now we start to prove the theorem. We divide the cases by $n$ and $a$.\\

\noindent \underline{\textit{Case 1.}}\quad Assume that $a=1$ or
$n=1$ or $n=a=2$. Then $(\P^n,\mathcal{O}_{\P^n}(a))$ satisfies
Property $N_p$ for all $p \geq 0$. For these cases, see Theorem
\ref{thm:1}, Theorem \ref{thm:2} and Theorem \ref{thm:3}.\\

\noindent \underline{\textit{Case 2.}}\quad Assume that $n \geq 3$
and $a=2$. Then $(\P^n,\mathcal{O}_{\P^n}(2))$ satisfies Property
$N_p$ if and only if $p \leq 5$. Thus it suffices to show that
there is an integer $s$ such that $(X,2H+\pi^*B)$ satisfies
Property $N_5$ for any $B \in \mbox{Pic}C$ such that
$\mbox{deg}(B) \geq s$. By Lemma \ref{lem:veryampleness}.$(3)$,
$\mbox{codim} (X,\P H^0 (X,L)) \geq 5$. Therefore $(X,2H+\pi^*B)$
satisfies Property $N_5$ if and only if
\begin{equation*}
H^1 (X,\wedge^{6} M_L \otimes L^j )=0~~\mbox{for all $j\geq1$.}
\end{equation*}
by Lemma \ref{lem:criterion}. First we show that $H^1
(X,\wedge^{6} M_L \otimes L^j )=0$ if $j \geq 3$. Indeed from the
short exact sequence
\begin{equation*}
0 \rightarrow \wedge^{7} M_L \otimes L^{j} \rightarrow \wedge^{7}
H^0 (X,L) \otimes L^{j} \rightarrow \wedge^6 M_L \otimes L^{j+1}
\rightarrow 0,
\end{equation*}
we obtain that $H^1 (X,\wedge^6 M_L \otimes L^{j+1}) \cong H^2 (X,
\wedge^7 M_L \otimes L^{j})$ for all $j \geq 1$ since $H^1
(X,L^j)=0$ for all $j \geq 1$. Note that $H^2 (X, \wedge^7 M_L
\otimes L^{j}) \subset H^2 (X, T^7 M_L \otimes L^{j})$ since
$\wedge^7 M_L$ is a direct summand of $T^7 M_L$. Also since $M_L$
is $1$ $\pi$-regular, $H^2 (X, T^7 M_L \otimes L^{j}) \cong H^2
(C,\pi_* T^7 M_L \otimes L^{j})=0$ for all $j \geq 3$ by Lemma
\ref{lem:directimage}. Therefore $(X,2H+\pi^*B)$ satisfies
Property $N_5$ if and only if
\begin{equation*}
H^1 (X,\wedge^{6} M_L \otimes L)=H^1 (X,\wedge^{6} M_L \otimes
L^2)=0.
\end{equation*}
Also from the short exact sequence
\begin{equation*}
0 \rightarrow \pi^* M_\mathcal{F} \rightarrow M_L \rightarrow
\mathcal{G}_2 \otimes \pi^* B  \rightarrow 0,
\end{equation*}
$H^1 (X,\wedge^6 M_L \otimes L^j)=0$ if
\begin{equation*}
H^1 (X,\wedge^i \pi^* M_{\mathcal{F}} \otimes \wedge^{6-i}
\mathcal{G}_2 \otimes \pi^* B^{6-i} \otimes L^j)=0~~\mbox{for all
$i=0,\cdots,6$.}
\end{equation*}
Now by $(*)$ and the projection formula,
\begin{equation*}
H^1 (X,\wedge^i \pi^* M_{\mathcal{F}} \otimes \wedge^{6-i}
\mathcal{G}_2 \otimes \pi^* B^{6-i} \otimes L^j) \cong H^1
(C,\wedge^i M_{\mathcal{F}} \otimes B^{6-i+j} \otimes \pi_*
(\wedge^{6-i} \mathcal{G}_2  \otimes 2jH)).
\end{equation*}
Thus by Lemma \ref{lem:folklore}.$(6)$ and Theorem
\ref{thm:Butlerestimation}, we get the desired vanishing if
\begin{eqnarray*}
& & \mu^- (\wedge^i M_{\mathcal{F}} \otimes B^{6-i+j}
\otimes \pi_* (\wedge^{6-i} \mathcal{G}_2  \otimes 2jH)) \\
& \geq & i \mu^- (M_{\mathcal{F}}) + (6-i+j) \mbox{deg}(B) + \mu^-
(\pi_*
(\wedge^{6-i} \mathcal{G}_2  \otimes 2jH)) \\
&\geq & - \frac{6\mu^- (\mathcal{F})}{\mu^- (\mathcal{F})-g}  +
(6-i+j)\mbox{deg}(B) + \mu^-
(\pi_*(\wedge^{6-i} \mathcal{G}_2  \otimes 2jH)) \\
& > & 2g-2.
\end{eqnarray*}
Since
\begin{equation*}
\mbox{min}  \{ \mu^- (\pi_*(\wedge^{6-i} \mathcal{G}_2  \otimes
2jH))~|~ i=0,\cdots,6~~\mbox{and}~~j=1,2 \}
\end{equation*}
is well-defined, there exists an integer $s$ such that
\begin{equation*}
- \frac{6\mu^- (\mathcal{F})}{\mu^- (\mathcal{F})-g}  +
(6-i+j)\mbox{deg}(B) + \mu^- (\pi_*(\wedge^{6-i} \mathcal{G}_2
\otimes 2jH))
> 2g-2
\end{equation*}
if $\mbox{deg}(B) \geq s$.\\

\noindent \underline{\textit{Case 3.}}\quad Assume that $n \geq 2$
and $a \geq 3$. In this case, the proof is almost the same as
\textit{Case 2}. Recall that if $(\P^n,\mathcal{O}_{\P^n}(a))$
satisfies Property $N_p$, then $p \leq 3a-3$\cite{OP}. By Lemma
\ref{lem:veryampleness}.$(3)$, $\mbox{codim} (X,\P H^0 (X,L)) \geq
3a-3$. Thus $(X,L)$ satisfies Property $N_p$ if and only if
\begin{equation*}
H^1 (X,\wedge^{p+1} M_L \otimes L^j )=0~~\mbox{for all $j\geq1$.}
\end{equation*}
But this can be proved for $j \geq 3$ by the same method as in
\textit{Case 2} since $p \leq 3a-3$. Also
\begin{equation*}
\mbox{min}  \{ \mu^- (\pi_*(\wedge^{p+1-i} \mathcal{G}_a  \otimes
ajH))~|~   0 \leq i \leq p+1~~\mbox{and}~~j=1,2 \}
\end{equation*}
is well-defined and hence there exists an integer $s(a,p)$ such
that
\begin{equation*}
H^1 (X,\wedge^{p+1} M_L \otimes L)=H^1 (X,\wedge^{p+1} M_L \otimes
L^2 )=0
\end{equation*}
if $\mbox{deg}(B) \geq s(a,p)$. \qed \\

Now we consider the converse of Theorem \ref{thm:positive}. One
may guess that if $(\P^n,\mathcal{O}_{\P^n}(a))$ fails to satisfy
Property $N_p$ then $(X,aH+\pi^*B)$ also fails to satisfy Property
$N_p$ no matter how large $\mbox{deg}(B)$ is. Note that this is
true for $C \times \P^n$ (i.e. $\mathcal{E} = \oplus^{n+1}
\mathcal{O}_C)$ by Elena Rubei's result\cite{Ru}. Recently D.
Eisenbud, M. Green, K. Hulek and S. Popescu prove the following:

\begin{theorem}[D. Eisenbud, M. Green, K. Hulek and S. Popescu, Theorem 1.1 in
\cite{EGHP}]\label{thm:EGHP} For a projective variety $X$ and a
very ample line bundle $L \in \mbox{Pic}X$, if $X \subset \P H^0
(X,L)$ admits a $(p+2)$-secant $p$-plane, i.e., there exists a
linear subspace $\Lambda \subset \P H^0 (X,L)$ of dimension $\leq
p$ such that $X \cap \Lambda$ is finite and $\mbox{length}(X \cap
\Lambda) \geq p+2$, then $(X,L)$ fails to satisfy Property $N_p$.
\end{theorem}

\noindent By applying this result to Veronese embedding, they
reproved
that\\
\begin{enumerate}
\item[$(i)$] if $n \geq 3$, then $(\P^n,\mathcal{O}_{\P^n}(2))$ fails
to satisfy Property $N_6$ since it admits a $8$-secant $6$-plane,
and
\item[$(ii)$] if $n \geq 2$ and $a \geq 3$, then
$(\P^n,\mathcal{O}_{\P^n}(a))$ fails to satisfy Property
$N_{3a-2}$ since it admits a $3a$-secant $(3a-2)$-plane.\\
\end{enumerate}

\noindent For details, see Proposition 3.2 in \cite{EGHP}. This
result enables us to prove the following:

\begin{corollary}\label{cor:negative}
Let $L = aH+\pi^* B$ be a very ample line bundle on $X$. For $\pi_* L = \mathcal{F}$,
assume that $\mathcal{O}_{\P_C (\mathcal{F})} (1)$ is also very ample. Then \\
$(1)$ If $n \geq 3$ and $a=2$, then $L$ fails to satisfy Property
$N_6$.\\
$(2)$ If $n \geq 2$ and $a \geq 3$, then $L$ fails to satisfy
Property $N_{3a-2}$.
\end{corollary}

\begin{proof}
The assertion comes immediately from Theorem \ref{thm:EGHP} since
$aH + \pi^* B$ defines a fiberwise $a$-uple Veronese embedding.
Indeed each fibre of $\pi$ admits a multisecant linear space as
described above.
\end{proof}

\section{Property $N_p$ of Ruled Varieties over a Curve}
\noindent Throughout this section we use the
following notations:\\
\begin{enumerate}
\item[$\bullet$] $C$ : smooth projective curve of genus $g$
\item[$\bullet$] $\mathcal{E}$ : vector bundle of rank $n+1$ over
$C$
\item[$\bullet$] $X = \P_C (\mathcal{E})$ with
$\mathcal{O}_{\P_C (\mathcal{E})} (1)=H$ and projection morphism
$\pi : X \rightarrow C$\\
\end{enumerate}
First recall that $(\P^n,\mathcal{O}_{\P^n} (a))$ satisfies
Property $N_p$ for all $p \geq 0$ if and only if $a=1$ or $n=1$
and $a \geq 1$ or $n=a=2$. For these three cases we obtain Theorem
\ref{thm:1}, Theorem \ref{thm:2} and Theorem \ref{thm:3}. Roughly
speaking, we prove that when $a=1$ or $n=1$ and $a \geq 1$ or
$n=a=2$ and if $p \geq 0$ is a fixed integer, then $(X,aH+\pi^*
B)$ satisfies Property $N_p$ for all $B \in \mbox{Pic}C$ having
sufficiently large degree. Note that this completes the proof of
Theorem \ref{thm:positive}. Also if $n \geq 3$ and $a=2$ or $n
\geq 2$ and $a \geq 3$, then $(\P^n,\mathcal{O}_{\P^n} (a))$
satisfies Property $N_p$ for $0 \leq p \leq a-1$. In these cases,
we prove an effective statement in Theorem \ref{thm:4}. Let
$L=aH+\pi^*B \in \mbox{Pic}X$ be such that $a \geq 1$ and $\mu^-
(\pi_* L) > 2g$. Thus $L$ is very ample and  $H^1 (X,L^j)=0$ for
all $j\geq 1$ by Lemma \ref{lem:veryampleness}. If we put $\pi_* L
= \mathcal{F}$, then by Lemma \ref{lem:directimage}.$(3)$ there is
an exact sequence
\begin{equation*}
0 \rightarrow  \mathcal{K}_{L} \rightarrow \pi^* \mathcal{F}
\rightarrow L \rightarrow 0
\end{equation*}
of vector bundles on $X$ where $\mathcal{K}_{L}$ is $1$
$\pi$-regular. So we have the following commutative
diagram:\\
\begin{equation}
\begin{CD}
&  &&  && 0 &\\
&  &&  && \downarrow &\\
& 0 && && \mathcal{K}_{L} &\\
& \downarrow &&  && \downarrow &\\
0 \rightarrow &\pi^* M_\mathcal{F}& \rightarrow & H^0 (C,\mathcal{F}) \otimes \mathcal{O}_X & \rightarrow  & \pi^* \mathcal{F} & \rightarrow 0 \\
& \downarrow && \parallel && \downarrow &\\
0 \rightarrow & M_L & \rightarrow & H^0 (X,L) \otimes \mathcal{O}_X & \rightarrow & L & \rightarrow 0 \\
& \downarrow &&  && \downarrow &\\
& \mathcal{K}_{L} && && 0  &\\
& \downarrow  &&  && &\\
& 0 &&  && &.\\
\end{CD}
\end{equation}\\

\begin{lemma}\label{lem:2}
If $\ell \leq aj$, then $\mu^- (\pi_* T^{\ell} \mathcal{K}_{L}
\otimes L^j) \geq (\ell+j)\mu^- (\mathcal{F})$.
\end{lemma}

\begin{proof}
Let $\mathcal{G}_a$ be the $1$ $\pi$-regular bundle on $X$ defined
by the short exact sequence
\begin{equation*}
0 \rightarrow \mathcal{G}_a \rightarrow \pi^* S^a (\mathcal{E})
\rightarrow aH \rightarrow 0.
\end{equation*}
Thus $\mathcal{K}_L = \mathcal{G}_a \otimes \pi^* B$. By applying
Lemma \ref{lem:directimage},
\begin{equation*}
\mu^- (\pi_* T^{\ell} \mathcal{G}_a \otimes ajH) \geq \ell \mu^-
(\pi_* \mathcal{G}_a \otimes H) + (aj-\ell) \mu^-(\mathcal{E}).
\end{equation*}
Also from two exact sequences
\begin{eqnarray*}
& & 0 \rightarrow \mathcal{G}_a \otimes H \rightarrow \pi^* \pi_*
(aH)\otimes H \rightarrow (a+1)H \rightarrow 0 \quad \mbox{and}\\
& & 0 \rightarrow \mathcal{G}_1 \otimes aH \rightarrow \pi^* \pi_*
(H) \otimes aH \rightarrow (a+1)H \rightarrow 0,
\end{eqnarray*}
$\pi_* (\mathcal{G}_a \otimes H)=\pi_* (\mathcal{G}_1 \otimes aH)$
because they are both the kernel of the map
\begin{equation*}
\mathcal{E} \otimes S^a (\mathcal{E}) \rightarrow S^{a+1}
\mathcal{E} \rightarrow 0.
\end{equation*}
Clearly $\mathcal{G}_1 = \Omega_{X/C} \otimes H$ and thus
\begin{equation*}
\mu^- (\pi_* \mathcal{G}_a \otimes H) = \mu^- (\pi_* \Omega_{X/C}
\otimes (a+1)H) \geq (a+1) \mu^- (\mathcal{E})
\end{equation*}
by Lemma \ref{lem:1}. Therefore $\mu^- (\pi_* T^{\ell}
\mathcal{G}_a \otimes ajH) \geq  a(\ell +j) \mu^-(\mathcal{E})$.
Note that the desired inequality is obtained immediately from this
since $\mathcal{K}_L = \mathcal{G}_a \otimes \pi^* B$.
\end{proof}

\subsection{Ruled Scrolls}
When $a=1$, we obtain the following:

\begin{theorem}\label{thm:1}
For $L=H+\pi^* B \in \mbox{Pic}X$, assume that $\mu^- (\pi_* L) >
2g$. Then $L$ satisfies Property $N_p$ if
\begin{equation*}
\mu^- (\pi_* L) + \frac{2g^2 -2g}{\mu^- (\pi_* L)}>  3g-1+p.
\end{equation*}
\end{theorem}

\begin{proof}
In the commutative diagram $(1)$, $\mathcal{F} = \mathcal{E}
\otimes B$ and $\mathcal{K}_L = \Omega_{X/C} \otimes L$ where
$\Omega_{X/C}$ denotes the relative canonical sheaf. Denote
$\wedge^j \Omega_{X/C}$ by $~\Omega_{X/C} ^j$. It needs to show
that $H^1 (X,\wedge^m M_L \otimes L^j)=0$ for $1 \leq m \leq p+1$
and all $j \geq 1$. Also from the sequence
\begin{equation*}
0 \rightarrow \pi^* M_\mathcal{F} \rightarrow M_L \rightarrow
\Omega_{X/C} \otimes L \rightarrow 0,
\end{equation*}
it suffices to show that for every $1 \leq \ell \leq p+1,~0 \leq i
\leq \ell$ and $j \geq 1$,
\begin{equation*}
H^1 (X,\wedge^i \pi^* M_{\mathcal{F}} \otimes \Omega_{X/C}
^{\ell-i} \otimes L^{\ell-i+j} )=0
\end{equation*}
or equivalently
\begin{equation*}
H^1 (C,\wedge^i  M_{\mathcal{F}} \otimes \pi_* (\Omega_{X/C}
^{\ell-i} \otimes L^{\ell-i+j}) )=0
\end{equation*}
by Proposition \ref{prop:Bott}. Thus by Lemma \ref{lem:folklore},
it suffices to show that
\begin{equation*}
\mu^- (\wedge^i M_{\mathcal{F}} \otimes \pi_* (\Omega_{X/C}
^{\ell-i} \otimes L^{\ell-i+j}) ) > 2g-2.
\end{equation*}
By the following Lemma \ref{lem:1},
\begin{equation*}
\mu^- (\pi_* \Omega_{X/C} ^{\ell-i} \otimes L^{\ell-i+j}) \geq
(\ell -i+j)\mu^- (\mathcal{F})
\end{equation*}
and hence
\begin{eqnarray*}
\mu^- (\wedge^i M_{\mathcal{F}} \otimes \pi_* (\Omega_{X/C}
^{\ell-i} \otimes L^{\ell-i+j}) )
& \geq & i  \mu^- (M_{\mathcal{F}}) + \mu^- (\pi_* (\Omega_{X/C}^{\ell-i} \otimes L^{\ell-i+j}))\\
 & \geq & i  \mu^- (M_{\mathcal{F}}) + (\ell -i+j)\mu^- (\mathcal{F}) \\
 & \geq & (p+1) \mu^- (M_{\mathcal{F}}) +  \mu^- (\mathcal{F})\\
 & \geq &  - \frac{(p+1)\mu^- (\mathcal{F})}{(\mu^- (\mathcal{F})-g)  }+  \mu^- (\mathcal{F})~~~~\mbox{(Theorem
 \ref{thm:Butlerestimation})}.
\end{eqnarray*}
It is easily checked that the inequality
\begin{equation*}
- \frac{(p+1)\mu^- (\mathcal{F})}{\mu^- (\mathcal{F})-g}+ \mu^-
(\mathcal{F}) > 2g-2
\end{equation*}
is equivalent to our assumption on $p$.
\end{proof}

\begin{lemma}\label{lem:1}
$\mu^- (\pi_* \Omega_{X/C} ^i \otimes L^{i+j}) \geq (i+j)\mu^-
(\mathcal{F})$ for every $j\geq1$.
\end{lemma}

\begin{proof}
In the commutative diagram $(1)$, $\mathcal{K}_L = \Omega_{X/C}
\otimes L$ and thus we have the short exact sequence
\begin{equation*}
0 \rightarrow \Omega_{X/C} \otimes L \rightarrow \pi^* \mathcal{F}
\rightarrow L \rightarrow 0.
\end{equation*}
From this we get the short exact sequence
\begin{equation*}
0 \rightarrow \Omega_{X/C}^{i+1} \otimes L^{i+j} \rightarrow
\wedge^{i+1} \pi^* \mathcal{F} \otimes L^{j-1} \rightarrow
\Omega_{X/C}^i \otimes L^{i+j} \rightarrow 0
\end{equation*}
for every $j \geq 1$. Then by Proposition \ref{prop:Bott}, $R^1
\pi_* (\Omega_{X/C}^{i+1} \otimes L^{i+j}) =0$ and hence we get
the short exact sequence
\begin{equation*}
0 \rightarrow \pi_*  \Omega_{X/C}^{i+1} \otimes L^{i+j}
\rightarrow \wedge^{i+1} \mathcal{F} \otimes \pi_* L^{j-1}
\rightarrow \pi_*  \Omega_{X/C}^i \otimes L^{i+j} \rightarrow 0.
\end{equation*}
Now by applying Lemma \ref{lem:folklore}.$(4)$,
\begin{equation*}
\mu^- (\pi_*  \Omega_{X/C}^i \otimes L^{i+j}) \geq \mu^- (
\wedge^{i+1} \mathcal{F} \otimes \pi_* L^{j-1}) \geq (i+j) \mu^-
(\mathcal{F}).
\end{equation*}
which completes the proof.
\end{proof}

\begin{corollary}\label{cor:detailscroll}
Let $L =H +\pi^* B \in \mbox{Pic}X$ be such that $\mbox{deg} (B)=b$.\\
$(1)$ When $g=1$, $L$ satisfies Property $N_p$ if $b+\mu^- (\mathcal{E}) > 2+p$.\\
$(2)$ When $g = 2$, $L$ satisfies Property $N_p$ if $b+\mu^- (\mathcal{E}) \geq 5+p$.\\
$(3)$ When $g =3$ and $0 \leq p \leq 4$, $L$ satisfies Property
$N_p$ if $b+\mu^- (\mathcal{E}) \geq 7+p$. \\
$(4)$ When $g =4$ and $0 \leq p \leq 2$, $L$ satisfies Property
$N_p$ if $b+\mu^- (\mathcal{E}) \geq 9+p$. \\
$(5)$ When $g =5$ and $0 \leq p \leq 2$, $L$ satisfies Property
$N_p$ if $b+\mu^- (\mathcal{E}) \geq 11+p$. \\
$(6)$ When $g \geq 2$, $L$ satisfies Property $N_p$ if $b+\mu^-
(\mathcal{E}) \geq 3g-1+p$.\\
$(7)$ When $g \geq 1$, $L$ is normally generated if $b+\mu^-
(\mathcal{E}) \geq 2g+1$ and normally

presented if $b+\mu^- (\mathcal{E}) \geq 2g+2$.
\end{corollary}

\begin{proof}
Put $\nu = \mu^- (\pi_* L) =  b+\mu^- (\mathcal{E})$. Then the
inequality in Theorem \ref{thm:1} is equivalent to
\begin{equation*}
\nu^2 -(3g-1+p)\nu +2g^2-2g > 0.
\end{equation*}
When $g=1$, this is equivalent to $\nu > 2+p$. When $g \geq 2$,
this holds for all $\nu \geq 3g-1+p$. The other cases are also
easily checked by using the above inequality.
\end{proof}

\noindent {\bf Remark 4.1.1.} If $g=0$ and $a=1$, i.e., $(X,L)$ is
a rational normal scroll, then Property $N_p$ holds for all $p
\geq 0$. See Theorem 3.6 in \cite{OP}. \qed \\

\noindent {\bf Example 4.1.1.} For $A_1,\cdots,A_n \in
\mbox{Pic}C$ of degrees $a_1,\cdots,a_n$ with $0 \geq a_1 \geq a_2
\geq \cdots \geq a_n$, let $\mathcal{E}= \mathcal{O}_C \oplus
\mathcal{O}_C (A_1) \oplus \cdots \oplus \mathcal{O}_C (A_n)$,
i.e., $\mathcal{E}$ is normalized and totally decomposable, and
let $X=\P_C (\mathcal{E})$. Thus $\mu^- (\mathcal{E})=a_n$. For $B
\in \mbox{Pic}C$ of degree $b$, $L=H+\pi^*B$ is very ample if
$b+a_n \geq 2g+1$. Also for the section $C_n$ determined by
$\mathcal{E} \rightarrow \mathcal{O}_C (A_n) \rightarrow 0$, the
restriction of $L$ to $C_n$ is equal to $\mathcal{O}_C (A_n + B)$.
Note that for the embedding $X \hookrightarrow \P H^0 (X,L)$,
$C_n$ is the linear section of $X$ and $\P H^0 (C,\mathcal{O}_C
(A_n + B)) \subset \P H^0 (X,L)$. For $g=1$ or $2$, we obtain the
followings:
\begin{enumerate}
\item[$(1)$] Let $C$ be an elliptic curve. Note that $L$ is very
ample if and only if $b+a_n \geq 3$. By Corollary
\ref{cor:detailscroll}, $(X,L)$ satisfies Property $N_p$ if $b+a_n
\geq 3+p$. Assume that $b+a_n = 3+p$. Then it fails to satisfy
Property $N_{p+1}$ since the linear section $C_n \subset \P H^0
(C,\mathcal{O}_C (A_n + B))$ gives a $(p+3)$-secant $(p+1)$-plane.
See Theorem \ref{thm:EGHP}. \item[$(2)$] Let $C$ be a curve of
genus $2$. Note that $L$ is very ample if and only if $b+a_n \geq
5$. By the same way as in $(1)$, $(X,L)$ satisfies Property $N_p$
if and only if $b+a_n \geq 5+p$.
\end{enumerate}
This seems to have been unknown even for the trivial case $X = C
\times \P^n$, i.e., $\mathcal{E} = \oplus ^{n+1}
\mathcal{O}_C$.\qed \\

\subsection{Ruled Surfaces}
Let $\mathcal{E}$ be a normalized rank $2$ vector bundle on $C$.
We follow the notation and terminology of R. Hartshorne's book
\cite{H}, V $\S 2$. Let
\begin{equation*}
\mathfrak{e}=\wedge^2 \mathcal{E}~~~~\mbox{and}~~~~e = -
\mbox{deg}(\mathfrak{e}).
\end{equation*}
We fix a minimal section $C_0$ such that $\mathcal{O}_X
(C_0)=\mathcal{O}_{\P_C (\mathcal{E})} (1)$. For $\mathfrak{b} \in
\mbox{Pic}C$, $\mathfrak{b}f$ denote the pullback of
$\mathfrak{b}$ by $\pi$. Thus any element of $\mbox{Pic}X$ can be
written $aC_0+\mathfrak{b}f$ with $a\in \Z$ and $\mathfrak{b} \in
\mbox{Pic}C$ and any element of $\mbox{Num}X$ can be written $aC_0
+bf$ with $a,b \in \Z$.

\begin{theorem}\label{thm:2}
For $L \equiv aC_0 + bf$, assume that $\mu^- (\pi_* L) > 2g$. Then
$L$ satisfies Property $N_p$ if
\begin{equation*}
a\mu^- (\mathcal{E})+b + \frac{2g^2 -2g}{a\mu^- (\mathcal{E})+b}>
3g-1 +p.
\end{equation*}
\end{theorem}

\begin{proof}
Note that in the commutative diagram $(1)$, $\mathcal{K}_{L}$ is a
vector bundle of rank $a$ and $1$ $\pi$-regular on $X$. It needs
to show that $H^1 (X,\wedge^m M_L \otimes L^j)=0$ for $1 \leq m
\leq p+1$ and all $j \geq 1$. Also from the sequence
\begin{equation*}
0 \rightarrow \pi^* M_\mathcal{F} \rightarrow M_L \rightarrow
\mathcal{K}_{L} \rightarrow 0,
\end{equation*}
it suffices to show that for every $1 \leq \ell \leq p+1,~0 \leq i
\leq \ell$ and $j \geq 1$,
\begin{equation*}
H^1 (X,\wedge^i \pi^* M_{\mathcal{F}} \otimes \wedge^{\ell-i}
\mathcal{K}_{L}  \otimes L^{j} )=0.
\end{equation*}\\
\noindent \underline{\textit{Case 1.}}\quad First we concentrate
on the case $\ell-i =0$. Then
\begin{equation*}
H^1 (X,\wedge^{\ell} \pi^*M_{\mathcal{F}} \otimes L^{j} )\cong H^1
(C,\wedge^{\ell} M_{\mathcal{F}} \otimes \pi_* L^{j})
\end{equation*}
and hence it suffices to show
\begin{equation*}
\mu^- (\wedge^{\ell} M_{\mathcal{F}} \otimes \pi_* L^{j}) > 2g-2.
\end{equation*}
Indeed we have
\begin{eqnarray*}
\mu^- (\wedge^{\ell} M_{\mathcal{F}} \otimes \pi_* L^{j}) & \geq &
\ell \mu^- (M_{\mathcal{F}})
+ j \mu^- ( \mathcal{F})\\
& \geq & (p+1) \mu^- (M_{\mathcal{F}})
+  \mu^- ( \mathcal{F})\\
& \geq & - (p+1) \frac{\mu^- (\mathcal{F})}{\mu^-
(\mathcal{F})-g}  + \mu^- (\mathcal{F}) \\
& > & 2g-2
\end{eqnarray*}
where the last inequality is equivalent to our assumption on
$p$ since
\begin{equation*}
\mu^- (\mathcal{F}) = a\mu^-(\mathcal{E})+b.
\end{equation*}

\noindent \underline{\textit{Case 2.}}\quad Now we consider the
case $\ell-i \geq 1$. Because $\mathcal{K}_L$ is of rank $a$, we
only consider the cases $\ell-i \leq a$. Since $\wedge^{\ell-i}
\mathcal{K}_{L}$ is a direct summand of the tensor product
$T^{\ell-i} \mathcal{K}_{L}$, it suffices to show that
\begin{equation*}
H^1 (X,\wedge^i \pi^* M_{\mathcal{F}} \otimes T^{\ell-i}
\mathcal{K}_{L} \otimes L^{j} )=0.
\end{equation*}
Note that since $\ell-i \leq a$, $T^{\ell-i} \mathcal{K}_{L}
\otimes L^{j}$ is $0$ $\pi$-regular for all $j\geq 1$. Therefore
\begin{equation*}
H^1 (X,\wedge^i \pi^* M_{\mathcal{F}} \otimes T^{\ell-i}
\mathcal{K}_{L} \otimes L^{j} )=H^1 (C,\wedge^i M_{\mathcal{F}}
\otimes \pi_* (T^{\ell-i} \mathcal{K}_{L} \otimes L^{j})).
\end{equation*}
By Lemma \ref{lem:2}, $\mu^- (\pi_* T^{\ell-i} \mathcal{K}_{L}
\otimes L^{j}) \geq (j+\ell -i) \mu^- (\mathcal{F})$. Therefore
\begin{eqnarray*}
\mu^- (\wedge^i M_{\mathcal{F}} \otimes \pi_* (T^{\ell-i}
\mathcal{K}_{L} \otimes L^{j}))
& \geq & (p+1) \mu^- (\mathcal{M}_{\mathcal{F}}) + \mu^- (\mathcal{F})\\
& > & 2g-2
\end{eqnarray*}
and hence we get the desired vanishing.
\end{proof}

\begin{corollary}\label{cor:detailsurface}
Let $L \in \mbox{Pic}X$ be a line bundle in the numerical class
$aC_0 +bf$.\\
$(1)$ When $g=1$, $L$ satisfies Property $N_p$ if $b+a \mu^-
(\mathcal{E}) > 2+p$.\\
$(2)$ When $g = 2$, $L$ satisfies Property $N_p$ if $b+a \mu^-
(\mathcal{E}) \geq 5+p$.\\
$(3)$ When $g =3$ and $0 \leq p \leq 4$, $L$ satisfies Property
$N_p$ if $b+a \mu^- (\mathcal{E}) \geq 7+p$. \\
$(4)$ When $g =4$ and $0 \leq p \leq 2$, $L$ satisfies Property
$N_p$ if $b+a \mu^- (\mathcal{E}) \geq 9+p$. \\
$(5)$ When $g =5$ and $0 \leq p \leq 2$, $L$ satisfies Property
$N_p$ if $b+a \mu^- (\mathcal{E}) \geq 11+p$. \\
$(6)$ When $g \geq 2$, $L$ satisfies Property $N_p$ if $b+a \mu^-
(\mathcal{E}) \geq 3g-1+p$.\\
$(7)$ When $g \geq 1$, $L$ is normally generated if $b+a \mu^-
(\mathcal{E}) \geq 2g+1$ and normally

presented if $b+a \mu^- (\mathcal{E}) \geq 2g+2$.
\end{corollary}

\begin{proof}
The results follow from the inequality of Theorem \ref{thm:2}. See
the proof of Corollary \ref{cor:detailscroll}.
\end{proof}

\noindent {\bf Remark 4.2.1.} Let X be the rational ruled surface
associated to $\mathcal{O}_{\P^1} \oplus \mathcal{O}_{\P^1} (-e)$.
For a given $p \geq 0$, F. J. Gallego and B. P.
Purnaprajna\cite{GP} classify all very ample line bundles on $X$
satisfying Property $N_p$. More precisely, let $L$ be a very ample
line bundle in the numerical class of $aC_0 +bf$(i.e., $a\geq1$
and $b-ae \geq 1$).
\begin{enumerate}
\item[$(i)$] If $a=1$, then $L$ satisfies Property $N_p$ for all
$p \geq 0$. \item[$(ii)$] If $e=0$ and $b=1$, then $L$ satisfies
Property $N_p$ for all $p \geq 0$. \item[$(iii)$] If $e=0$ and
$a,b \geq 2$ or $e \geq 1$ and $a \geq 2$, then $L$ satisfies
Property $N_p$ if and only if $a \geq 1$ and $2a+2b-ae \geq 3+p$.\qed \\
\end{enumerate}

\noindent {\bf Remark 4.2.2.} In \cite{Park}, the relation between
higher syzygies of the ruled surface and that of the minimal
section is considered. Also it is proved that Mukai's
conjecture is true for ruled surfaces with $e \geq \mbox{max}\{0,\frac{g-2}{2}\}$. \qed \\

\subsection{Veronese Surface Fibrations}
Let $\mathcal{E}$ is a rank $3$ vector bundle and as in the
previous subsection, we use the following notations:
\begin{equation*}
\mathfrak{e} = \wedge^3 \mathcal{E} ~~~~\mbox{and}~~~~e = -
\deg(\mathfrak{e} )
\end{equation*}

\begin{theorem}\label{thm:3}
For $L=2H+\pi^* B \in \mbox{Pic}X$, assume that $\mu^- (\pi_* L) >
2g$. Then $L$ satisfies Property $N_p$ if
\begin{equation*}
7 \mu (\pi_* L) \geq \mu^+ (\pi_* L)~~~\mbox{and}~~~~\mu^- (\pi_*
L) + \frac{2g^2 -2g}{\mu^- (\pi_* L)} >  3g-1+p.
\end{equation*}
$($Note that when $\mathcal{E}$ is semistable, the first condition
always holds .$)$
\end{theorem}

\begin{proof}
Again we use the commutative diagram $(1)$ introduced in the
first of this section. From the exact sequence
\begin{equation*}
0 \rightarrow \pi^* M_{\mathcal{F}} \rightarrow \mathcal{M}_{L}
\rightarrow \mathcal{K}_{L} \rightarrow 0,
\end{equation*}
it suffices to show that
\begin{equation*}
H^1 (X,\wedge^i \pi^* M_{\mathcal{F}} \otimes \wedge^{\ell-i}
\mathcal{K}_{L} \otimes L^{j} )=0~~~~\mbox{for every $1 \leq \ell
\leq p+1,~0 \leq i \leq \ell$ and $j \geq 1$.}
\end{equation*}
But since $\mathcal{K}_{L}$ is a rank $5$ vector bundle on $X$, we
only consider the cases $\ell-i \leq 5$.\\

\noindent \underline{\textit{Case 1.}}\quad Assume that $\ell-i
=0$. We want to prove that
\begin{equation*}
H^1 (X,\wedge^i \pi^* M_{\mathcal{F}} \otimes L^{j}
)=0~~~~\mbox{for every $1 \leq i \leq p+1$ and $j \geq 1$}.
\end{equation*}
But $H^1 (X,\wedge^i \pi^* M_{\mathcal{F}} \otimes L^{j} ) \cong
H^1 (C,\wedge^i M_{\mathcal{F}} \otimes \pi_* L^{j} )$ and
\begin{eqnarray*}
\mu^- (\wedge^i M_{\mathcal{F}} \otimes \pi_* L^{j}) & \geq &
(p+1) \mu^- (M_{\mathcal{F}})
+ \mu^- (\mathcal{F}) \\
& \geq & -(p+1) \frac{\mu^- (\mathcal{F})}{(\mu^-
(\mathcal{F})-g)}  + \mu^- (\mathcal{F}) \\
& > &  2g-2
\end{eqnarray*}
which gives the desired vanishing.\\

\noindent \underline{\textit{Case 2.}}\quad Assume that $1 \leq
\ell-i \leq 2$. Here we check that
\begin{equation*}
H^1 (X,\wedge^i \pi^* M_{\mathcal{F}} \otimes \wedge^{k}
\mathcal{K}_{L} \otimes L^{j} )=0~~~~\mbox{for every $1\leq
k\leq2$, $1 \leq i+k \leq p+1$ and $j \geq 1$}.
\end{equation*}
Also since $\wedge^{k} \mathcal{K}_{L}$ is a direct summand of the
tensor product $T^{k} \mathcal{K}_{L}$, we may instead show that
\begin{equation*}
H^1 (X,\wedge^i \pi^* M_{\mathcal{F}} \otimes T^{k}
\mathcal{K}_{L} \otimes L^{j} )=0.
\end{equation*}
In this case, $T^{k}\mathcal{K}_{L}\otimes L^{j}$ is $0$
$\pi$-regular for $j\geq 1$ and hence the above vanishing is
obtained by a similar method as in the proof of
Theorem \ref{thm:2}.\\

\noindent \underline{\textit{Case 3.}}\quad Assume that $\ell-i =
3$. What we want to show is
\begin{equation*}
H^1 (X,\wedge^i \pi^* M_{\mathcal{F}} \otimes \wedge^{3}
\mathcal{K}_{L} \otimes L^{j} )=0~~~~\mbox{for every $1 \leq i
\leq p-2$ and $j \geq 1$}.
\end{equation*}
Note that if $j\geq 2$, then $\wedge^{3} \mathcal{K}_{L} \otimes
L^{j}$ is $0$ $\pi$-regular and we can use the same method as in
the second case. Also $\wedge^{3} \mathcal{K}_{L} \otimes L$ is
$1$ $\pi$-regular and hence we concentrate on proving the
following:
\begin{equation*}
H^1 (C,\wedge^i M_{\mathcal{F}} \otimes \pi_* \wedge^{3}
\mathcal{K}_{L} \otimes L)=0
\end{equation*}
But by the following Lemma \ref{lem:complicated1}, $\mu^- (\pi_*
\wedge^{3} \mathcal{K}_{L} \otimes L) \geq \mu^- (\mathcal{F})$
and hence
\begin{equation*}
\mu^- (\wedge^i M_{\mathcal{F}} \otimes \pi_* \wedge^{3}
\mathcal{K}_{L} \otimes L) \geq (p+1) \mu^- (M_{\mathcal{F}}) +
\mu^- (\mathcal{F}) > 2g-2
\end{equation*}
which completes the proof.\\

\noindent \underline{\textit{Case 4.}}\quad Assume that $\ell-i =
4$. The following is our aim:
\begin{equation*}
H^1 (X,\wedge^i \pi^* M_{\mathcal{F}} \otimes \wedge^{4}
\mathcal{K}_{L} \otimes L^{j} )=0~~~~\mbox{for every $1 \leq i
\leq p-3$ and $j \geq 1$}.
\end{equation*}
Note that if $j\geq 2$, then $\wedge^{3} \mathcal{K}_{L} \otimes
L^{j}$ is $0$ $\pi$-regular and we can use the same method as in
the second case. Thus it suffices to show that
\begin{equation*}
H^1 (X,\wedge^i \pi^* M_{\mathcal{F}} \otimes \wedge^{4}
\mathcal{K}_{L} \otimes L)=0~~~~\mbox{for every $1 \leq i \leq
p-3$}.
\end{equation*}
From the short exact sequence $0 \rightarrow \mathcal{K}_{L}
\rightarrow \pi^* \mathcal{F} \rightarrow L \rightarrow 0$ we have
\begin{equation*}
0 \rightarrow \wedge^{5} \mathcal{K}_{L} \rightarrow \wedge^5
\pi^* \mathcal{F} \rightarrow \wedge^{4} \mathcal{K}_{L} \otimes L
\rightarrow 0
\end{equation*}
and $\wedge^{5} \mathcal{K}_{L} \cong \wedge^6 \pi^* \mathcal{F}
\otimes L^{-1}$. In particular, $R^1 \pi_* \wedge^{5}
\mathcal{K}_{L} = R^2 \pi_* \wedge^{5} \mathcal{K}_{L} =0$ since
$L=2H+\pi^* B$. Therefore
\begin{equation*}
H^1 (X,\wedge^i \pi^* M_{\mathcal{F}} \otimes \wedge^{4}
\mathcal{K}_{L} \otimes L) \cong H^1 (C,\wedge^i M_{\mathcal{F}}
\otimes \pi_* \wedge^{4} \mathcal{K}_{L} \otimes L).
\end{equation*}
Also from the short exact sequence
\begin{equation*}
0 \rightarrow \pi_* \wedge^{5} \mathcal{K}_{L} \rightarrow
\wedge^5  \mathcal{F} \rightarrow \pi_* \wedge^{4} \mathcal{K}_{L}
\otimes L \rightarrow 0,
\end{equation*}
$\mu^- (\pi_* \wedge^{4} \mathcal{K}_{L} \otimes L) \geq 5 \mu^-
(\mathcal{F})$ and hence
\begin{equation*}
\mu^- (\wedge^i M_{\mathcal{F}} \otimes \pi_* \wedge^{4}
\mathcal{K}_{L} \otimes L) \geq (p+1) \mu^- (M_{\mathcal{F}}) +
\mu^- (\mathcal{F}) > 2g-2
\end{equation*}
which completes the proof.\\

\noindent \underline{\textit{Case 5.}}\quad Assume that $\ell-i =
5$. Since $\wedge^5 \mathcal{K}_{L} = \wedge^6 \pi^* \mathcal{F}
\otimes L^{-1}$, we need to check that
\begin{equation*}
H^1 (X,\wedge^i \pi^* M_{\mathcal{F}} \otimes  \wedge^6 \pi^*
\mathcal{F}  \otimes L^{j-1} )=0~~~~\mbox{for every $1 \leq i
\leq p-4$ and $j \geq 1$}.
\end{equation*}
Since $H^1 (X,\wedge^i \pi^* M_{\mathcal{F}} \otimes \wedge^6
\pi^* \mathcal{F}  \otimes L^{j-1} ) \cong H^1 (C,\wedge^i
M_{\mathcal{F}} \otimes  \wedge^6 \mathcal{F}  \otimes \pi_*
L^{j-1} )$, the desired vanishing follows immediately by
estimating $\mu^-$.
\end{proof}

\begin{lemma}\label{lem:complicated1}
$\mu^- (\pi_* \wedge^{3} \mathcal{K}_{L} \otimes L) \geq \mu^-
(\mathcal{F})$.
\end{lemma}

\begin{proof}
From the short exact sequence
\begin{equation*}
0 \rightarrow \wedge^5 \mathcal{K}_{L} \otimes L^{-1} \cong
\wedge^6 \pi^* \mathcal{F} \otimes L^{-2} \rightarrow \wedge^5
\pi^* \mathcal{F}\otimes L^{-1} \rightarrow \wedge^{4}
\mathcal{K}_{L} \rightarrow 0.
\end{equation*}
we have
\begin{equation*}
\begin{CD}
0 = R^1 \pi_* \wedge^5 \pi^* \mathcal{F}\otimes L^{-1} &
\rightarrow & R^1 \pi_* \wedge^{4} \mathcal{K}_{L} \rightarrow R^2
\pi_* \wedge^6 \pi^* \mathcal{F} \otimes L^{-2} \\ & \rightarrow &
R^2 \pi_* \wedge^5 \pi^* \mathcal{F}\otimes L^{-1}=0.
\end{CD}
\end{equation*}
Therefore $R^1 \pi_* \wedge^{4} \mathcal{K}_{L} \cong R^2 \pi_*
\wedge^6 \pi^* \mathcal{F} \otimes L^{-2} \cong \mathcal{E}^{\vee}
\otimes (\wedge^3 \mathcal{E})^{\vee} \otimes \mathcal{O}_C (-2B)
\otimes \wedge^6 \mathcal{F}$. For details, see Exercise 3.8.4 in
\cite{H}. Also the short exact sequence $0 \rightarrow \wedge^{4}
\mathcal{K}_{L} \rightarrow \wedge^4 \pi^* \mathcal{F} \rightarrow
\wedge^{3} \mathcal{K}_{L} \otimes L \rightarrow 0$ induces
\begin{eqnarray*}
0 \rightarrow \pi_* \wedge^{4} \mathcal{K}_{L} & \rightarrow &
\pi_* \wedge^4 \pi^* \mathcal{F} \rightarrow  \pi_*
\wedge^{3} \mathcal{K}_{L} \otimes L \\
& \rightarrow & R^1 \pi_* \wedge^{4} \mathcal{K}_{L} \rightarrow
R^1 \wedge^4 \pi^* \mathcal{F} =0.
\end{eqnarray*}
Let the image of $\pi_* \wedge^4 \pi^* \mathcal{F} \rightarrow
\pi_* \wedge^{3} \mathcal{K}_{L} \otimes L$ be $\mathcal{H}$. Then
from the diagram
\begin{equation*}
\begin{CD}
& 0 & \\
& \downarrow & \\
0 \rightarrow \pi_* \wedge^{4} \mathcal{K}_{L} \rightarrow  \pi_* \wedge^4 \pi^* \mathcal{F} \rightarrow & \mathcal{H} & \rightarrow 0 \\
& \downarrow & \\
& \pi_* \wedge^{3} \mathcal{K}_{L} \otimes L & \\
& \downarrow & \\
& \mathcal{E}^{\vee} \otimes (\wedge^3 \mathcal{E})^{\vee} \otimes
(-2B) \otimes \wedge^6 \mathcal{F} &\\
& \downarrow & \\
& 0 & \\
\end{CD}
\end{equation*}
we have
\begin{eqnarray*}
\mu^- (\pi_* \wedge^{3} \mathcal{K}_{L} \otimes L) & \geq &
\mbox{min}\{ \mu^- (\mathcal{H}),\mu^- ( \mathcal{E}^{\vee}
\otimes (\wedge^3 \mathcal{E})^{\vee} \otimes
(-2B) \otimes \wedge^6 \mathcal{F}) \} \\
& \geq & \mbox{min}\{ 4\mu^- (\mathcal{F}),-\mu^+  ( \mathcal{E})
+e  -2 \deg(B) + 6 \mu (\mathcal{F})) \}\\
& = & \mbox{min}\{ 4\mu^- (\mathcal{F}),-\mu^+  ( \mathcal{E})
  - \frac{1}{2} \deg(B) +  \frac{9}{2}  \mu (\mathcal{F})) \}\\
& = & \mbox{min}\{ 4\mu^- (\mathcal{F}),-\frac{1}{2}\mu^+  (
\mathcal{F}) +  \frac{9}{2}  \mu (\mathcal{F})) \} \\
& \geq &  \mu^- (\mathcal{F}).
\end{eqnarray*}
Here we use the identities $\mu^- ( \mathcal{E}^{\vee} )=-\mu^+
(\mathcal{E})$, $\deg S^2 \mathcal{E} = 4\deg \mathcal{E}$ and $e
= \frac{3}{2}\deg(B ) - \frac{3}{2} \mu (\mathcal{F})$, and the
assumption $7 \mu (\mathcal{F}) \geq \mu^+  ( \mathcal{F})$.
\end{proof}

\subsection{Arbitrary Cases} In this subsection we consider arbitrary
$n  \geq 1$ and $a \geq 1$ and $0 \leq p \leq a-1$. Recall that
$(\P^n,\mathcal{O}_{\P^n} (a))$ satisfies Property $N_p$ for $0
\leq p \leq a-1$. Precisely we prove the following:

\begin{theorem}\label{thm:4}
For $L=aH+\pi^* B \in \mbox{Pic}X$ , assume that $\mu^- (\pi_* L)
> 2g$. If $0 \leq p \leq a-1$ and
\begin{equation*}
\mu^- (\pi_* L) + \frac{2g^2 -2g}{\mu^- (\pi_* L)}>  3g-1+p,
\end{equation*}
then $(X,L)$ satisfies Property $N_p$.
\end{theorem}

\begin{proof}
Again we use the short exact sequence
\begin{equation*}
0 \rightarrow \pi^* M_{\mathcal{F}} \rightarrow \mathcal{M}_{L}
\rightarrow \mathcal{K}_{L} \rightarrow 0,
\end{equation*}
in the commutative diagram $(1)$. Thus it suffices to show that
\begin{equation*}
H^1 (X,\wedge^i \pi^* M_{\mathcal{F}} \otimes \wedge^{\ell-i}
\mathcal{K}_{L} \otimes L^{j} )=0~~~~\mbox{for every $1 \leq \ell
\leq p+1,~0 \leq i \leq \ell$ and $j \geq 1$.}
\end{equation*}
Since $\wedge^{\ell-i} \mathcal{K}_{L}$ is a direct summand of the
tensor product $T^{\ell-i} \mathcal{K}_{L}$, we prove that
\begin{equation*}
H^1 (X,\wedge^i \pi^* M_{\mathcal{F}} \otimes T^{\ell-i}
\mathcal{K}_{L} \otimes L^{j} )=0.
\end{equation*}
Note that since $\ell-i \leq p+1 \leq a$, $T^{\ell-i}
\mathcal{K}_{L} \otimes L^{j}$ is $0$ $\pi$-regular for all $j\geq
1$. Therefore
\begin{equation*}
H^1 (X,\wedge^i \pi^* M_{\mathcal{F}} \otimes T^{\ell-i}
\mathcal{K}_{L} \otimes L^{j} )=H^1 (C,\wedge^i M_{\mathcal{F}}
\otimes \pi_* (T^{\ell-i} \mathcal{K}_{L} \otimes L^{j})).
\end{equation*}
By Lemma \ref{lem:2}, $\mu^- (\pi_* (T^{\ell-i} \mathcal{K}_{L}
\otimes L^{j})) \geq (\ell -i +j) \mu^- (\mathcal{F})$ and hence
\begin{eqnarray*}
\mu^- (\wedge^i M_{\mathcal{F}} \otimes \pi_* (T^{\ell-i}
\mathcal{K}_{L} \otimes L^{j})) & \geq & (p+1)\mu^- (\wedge^i
M_{\mathcal{F}})+\mu^- (\mathcal{F}) > 2g-2
\end{eqnarray*}
by our assumption on $p$, which completes the proof.
\end{proof}

\begin{corollary}\label{cor:detailarbitrary}
Let $L =aH +\pi^* B \in \mbox{Pic}X$ be such that $\mbox{deg} (B)=b$.
For $(1)- (6)$, we always assume that $p \leq a-1$\\
$(1)$ When $g=1$ , $L$ satisfies Property $N_p$ if $b+a\mu^- (\mathcal{E}) > 2+p$.\\
$(2)$ When $g = 2$, $L$ satisfies Property $N_p$ if $b+a\mu^- (\mathcal{E}) \geq 5+p$.\\
$(3)$ When $g =3$ and $0 \leq p \leq 4$, $L$ satisfies Property
$N_p$ if $b+a\mu^- (\mathcal{E}) \geq 7+p$. \\
$(4)$ When $g =4$ and $0 \leq p \leq 2$, $L$ satisfies Property
$N_p$ if $b+a\mu^- (\mathcal{E}) \geq 9+p$. \\
$(5)$ When $g =5$ and $0 \leq p \leq 2$, $L$ satisfies Property
$N_p$ if $b+a\mu^- (\mathcal{E}) \geq 11+p$. \\
$(6)$ When $g \geq 2$, $L$ satisfies Property $N_p$ if $b+a\mu^-
(\mathcal{E}) \geq 3g-1+p$.\\
$(7)$ When $g \geq 1$, $L$ is normally generated if $b+a\mu^-
(\mathcal{E}) \geq 2g+1$ and normally

presented if $b+a\mu^- (\mathcal{E}) \geq 2g+2$.
\end{corollary}

\begin{proof}
The results follow from the inequality of Theorem \ref{thm:4}. See
the proof of Corollary \ref{cor:detailscroll}.
\end{proof}

\noindent {\bf Example 4.4.1.} We consider again the ruled variety
defined in Example 4.1.1. For $B \in \mbox{Pic}C$ of degree $b$,
$L=aH+\pi^*B$ is very ample if $b+a\cdot a_n \geq 2g+1$. Also for
the section $C_n$ determined by $\mathcal{E} \rightarrow
\mathcal{O}_C (A_n) \rightarrow 0$, the restriction of $L$ to
$C_n$ is equal to $\mathcal{O}_C (aA_n + B)$. Note that for the
embedding $X \hookrightarrow \P_C (\pi_* L) \hookrightarrow \P H^0
(X,L)$, $C_n$ is the linear section of $X$ and $\P H^0
(C,\mathcal{O}_C (aA_n + B)) \subset \P H^0 (X,L)$. For $g=1,2$,
we obtain the followings:
\begin{enumerate}
\item[$(1)$] Let $C$ be an elliptic curve. Note that $L$ is very
ample if and only if $b+a \cdot a_n \geq 3$. By Corollary
\ref{cor:detailarbitrary}, if $a \geq p+1$, then $(X,L)$ satisfies
Property $N_p$ for $b+a \cdot a_n \geq 3+p$. Assume that $b+a
\cdot a_n = 3+p$. Then it fails to satisfy Property $N_{p+1}$
since the linear section $C_n \subset \P H^0 (C,\mathcal{O}_C
(aA_n + B))$ gives a $(p+3)$-secant $(p+1)$-plane. See Theorem
\ref{thm:EGHP}. \item[$(2)$] Let $C$ be a curve of genus $2$. Note
that $L$ is very ample if and only if $b+a \cdot a_n \geq 5$. By
the same way as in $(1)$, if $a \geq p+1$, then $(X,L)$ satisfies
Property $N_p$ if and only if $b+a \cdot a_n \geq 5+p$.\qed \\
\end{enumerate}

\begin{corollary}\label{cor:Mukai}
For a vector bundle $\mathcal{E}$ of rank $n$ on $C$, put $\mu^-
(\mathcal{E}) = \frac{\nu}{\tau}$ where $\nu \in \Z$, $\tau \in
\N$ and $(\nu,\tau)=1$. Let $X = \P_C (\mathcal{E})$ be the
associated projective space bundle and let $L=K_X + A_1 + \cdots
+A_q$ where $A_i \in \mbox{Pic}X$ is ample and $q \geq n+1$. Then
\begin{equation*}
\mbox{$L$ satisfies Property $N_p$ if $q > \tau (g+1+p)$.}
\end{equation*}
\end{corollary}

\begin{proof}
Let $H = \mathcal{O}_{\P_C (\mathcal{E})} (1)$  and put $A_i = a_i
H+ \pi^* B_i$ where $a_i \geq 1$ and $B_i \in \mbox{Pic}C$. Then
\begin{equation*}
K_X + A_1 + \cdots +A_q = (\sum a_i -n) H + \pi^* (\sum B_i + K_C
+ \wedge^n \mathcal{E})
\end{equation*}
Therefore
\begin{equation*}
\mu^- (\pi_* L) \geq  \sum  \mu^- (\pi_* A_i) + 2g-2 + n
(\mu(\mathcal{E}) - \mu^- (\mathcal{E})).
\end{equation*}
Since $A_i$ is ample, $a_i \geq 1$ and $\mu^- (\pi_* A_i) >0$ by
Lemma \ref{lem:ampleMiyaoka}. Thus $\sum a_i -n \geq q-n \geq 1$.
Also since $\mu^- (\mathcal{E})=\frac{\nu}{\tau}$, we have $\mu^-
(\pi_* A_i)=a_i \mu^- (\mathcal{E}) + \mbox{deg}(B_i) \geq
\frac{1}{\tau}$. Thus if
\begin{equation*}
\mu^- (\pi_* L) \geq  2g-2 + \frac{q}{\tau} > 3g-1+p,
\end{equation*}
then $L$ satisfies Property $N_p$ by Theorem \ref{thm:4}. Clearly
the second inequality is equivalent to $q > \tau (g+1+p)$.
\end{proof}

\noindent {\bf Remark 4.4.1.} If $\mu^- (\mathcal{E})$ is an
integer(e.g. $\mathcal{E}$ is totally decomposable), then
Corollary \ref{cor:Mukai} says that $K_X + A_1 + \cdots +A_q$
satisfies Property $N_p$ if $q \geq \mbox{max} \{ n+1,g+2+p \}$.
In particular, Mukai's
conjecture is true for $X= \P_C (\mathcal{E})$ if $n \geq g$.\qed \\

\noindent {\bf Example 4.4.2.} Let $\mathcal{E}$ be a vector
bundle of rank $n$ on $\P^1$ and let $X = \P_{\P^1}
(\mathcal{E})$. Since every vector bundle on $\P^1$ is totally
decomposable, $K_X + A_1 + \cdots
+A_q$ satisfies Property $N_p$ if $q \geq \mbox{max} \{ n+1,2+p \}$. \qed \\

\section{Optimality of the condition of Mukai's conjecture}
\noindent This section is devoted to construct examples which show
that the condition of Mukai's conjecture is sharp. Indeed for any
$n \geq 2$ and $p \geq 0$ we construct $(X,A)$ where $X$ is a
smooth projective variety of dimension $n$ and $A \in \mbox{Pic}
X$ is an ample line bundle such that $K_X + (n+2+p)A$ fails to
satisfy Property $N_{p+1}$. Let $C$ be a smooth projective curve
of genus $g \geq 1$.

\begin{lemma}\label{lem:extension}
Let $\mathcal{E}$ be a normalized semistable vector bundle over
$C$ of rank $n \geq 2$ and $\mbox{deg} ( \wedge^n \mathcal{E})=1$.
For $\xi \in Ext^1 _{\mathcal{O}_C} (\mathcal{E},\mathcal{O}_C)$,
consider the corresponding extension
\begin{equation*}
 0 \rightarrow \mathcal{O}_C \rightarrow \mathcal{E}'
\rightarrow \mathcal{E} \rightarrow 0 \quad \cdots(\xi).
\end{equation*}
$(1)$ $\mathcal{E}'$ is normalized and $\mbox{deg} ( \wedge^n
\mathcal{E})=1$.\\
$(2)$ If $\xi \neq 0$, then $\mathcal{E}'$ is semistable.
\end{lemma}

\begin{proof}
$(1)$ Obvious.\\
$(2)$ If $\mathcal{E}'$ is unstable, there is a quotient $0
\rightarrow \mathcal{F} \rightarrow \mathcal{E}' \rightarrow
\mathcal{G} \rightarrow 0$ where $\mathcal{G}$ is a semistable
vector bundle and
\begin{equation*}
\frac{1}{n+1}= \mu (\mathcal{E}') > \mu^- (\mathcal{E}') = \mu
(\mathcal{G}) =
\frac{\mbox{deg}(\mathcal{G})}{\mbox{rank}(\mathcal{G})}.
\end{equation*}
Therefore $\mbox{deg}(\mathcal{G}) \leq 0$.\\

\noindent \underline{\textit{Step 1.}}\quad If
$\mbox{deg}(\mathcal{G})<0$, then $H^1 (C,\mathcal{G})=0$ since
$\mu^+ (\mathcal{G}) = \mu (\mathcal{G}) =\mu^- (\mathcal{G}) <0$.
Thus the composed map $\mathcal{O}_C \rightarrow \mathcal{G}$ is
the zero map. This implies that the surjection $\mathcal{E}'
\rightarrow \mathcal{G} \rightarrow 0$ factors through
$\mathcal{E} \rightarrow \mathcal{G} \rightarrow 0$. Since
$\mathcal{E}$ is semistable, $\mu (\mathcal{E})=\frac{1}{n} \leq
\mu (\mathcal{G}) < 0$ which is a contradiction. Therefore we
assume that
$\mbox{deg}(\mathcal{G})=0$ and $\mbox{deg}(\mathcal{F})=1$.\\

\noindent \underline{\textit{Step 2.}}\quad If the composed map
$\mathcal{O}_C \rightarrow \mathcal{G}$ is the zero map, then we
obtain a contradiction by the same way. Thus we assume that
$\mathcal{O}_C \rightarrow \mathcal{G}$ is injective. Now consider
the following commutative diagram:
\begin{equation*}
\begin{CD}
              &               &             &   0              &             &   0         &               \\
              &               &             &  \downarrow      &             &  \downarrow &               \\
              &               &             & \mathcal{F}      & =           & \mathcal{F} &      \\
              &               &             & \downarrow       &             & \downarrow  &               \\
0 \rightarrow & \mathcal{O}_C & \rightarrow &  \mathcal{E}'    & \rightarrow & \mathcal{E} & \rightarrow  0 \\
              &  \parallel    &             & \downarrow       &             & \downarrow  &               \\
0 \rightarrow & \mathcal{O}_C & \rightarrow &  \mathcal{G}     & \rightarrow & \mathcal{G}/\mathcal{O}_C & \rightarrow 0   \\
              &               &             &  \downarrow      &             &  \downarrow &               \\
              &               &             &  0               &             &  0              &
              \\\\
\end{CD}
\end{equation*}
Since $\mathcal{E}$ is semistable,
\begin{equation*}
\mu (\mathcal{F}) = \frac{1}{\mbox{rank} \mathcal{F}} \leq
\frac{1}{n} = \mu (\mathcal{E})
\end{equation*}
while $\mbox{rank}(\mathcal{F}) \leq n$. Therefore
$\mbox{rank}(\mathcal{F}) = n$ and $\mathcal{G}/\mathcal{O}_C$ is
a torsion sheaf. This implies that $\mathcal{G}$ is a line bundle
of degree $0$. Thus the injection $0 \rightarrow \mathcal{O}_C
\rightarrow \mathcal{G}$ is indeed an isomorphism, $\mathcal{F}
\cong \mathcal{E}$ and hence $\xi =0$ which is a contradiction.
\end{proof}

\noindent When $\mathcal{E}$ is normalized, $Ext^1
_{\mathcal{O}_C} (\mathcal{E},\mathcal{O}_C) \cong H^1
(C,\mathcal{E}^{\vee}) \cong H^0 (C,\mathcal{E} \otimes
K_C)^{\vee} \neq 0$. Also there exists a normalized semistable
rank $2$ vector bundle of degree $1$. That is, there is a ruled
surface over $C$ with $e=-1$. Thus by Lemma \ref{lem:extension},
for any $n \geq 2$ we have a sequence of vector bundles $\{
\mathcal{E}_2,\cdots,\mathcal{E}_n \}$ such that\\
\begin{enumerate}
\item[$(\alpha)$] $\mathcal{E}_i$ is a normalized semistable rank $i$
vector bundle of degree $1$ and
\item[$(\beta)$] there is a short
exact sequence $0 \rightarrow \mathcal{O}_C \rightarrow
\mathcal{E}_{i+1} \rightarrow \mathcal{E}_i \rightarrow 0$.\\
\end{enumerate}
Now let $X_i = \P_C (\mathcal{E}_i)$. Then by Lemma
\ref{lem:ampleMiyaoka}, $H_i = \mathcal{O}_{\P_C (\mathcal{E}_i)}
(1)$ is ample since $\mu^- (\mathcal{E}_i)=1/i
>0$. Also for the embedding $X_i \hookrightarrow X_{i+1}$ given by the sequence $0
\rightarrow \mathcal{O}_C \rightarrow \mathcal{E}_{i+1}
\rightarrow \mathcal{E}_i \rightarrow 0$, $H_{i+1} |_{X_i} = H_i$.

\begin{theorem}
Assume that $C$ is a hyperelliptic curve. Then for each $p \geq 0$,\\
$K_{X_n} + (n+2+p)H_n$ fails to satisfy Property $N_{p+1}$.
\end{theorem}

\begin{proof}
To simplify notations, put $L_{i,p} = K_{X_{i}}+(i+2+p)H_{i} \in
\mbox{Pic} (X_i)$ and
\begin{equation*}
0 \rightarrow M_{i,p} \rightarrow H^0 (X_i,L_{i,p}) \otimes
\mathcal{O}_{X_{i}} \rightarrow L_{i,p} \rightarrow 0.
\end{equation*}
Observe that $L_{i+1,p}|_{X_i} = L_{i,p}$. Clearly it suffices to
show the following
two statements:\\
\begin{enumerate}
\item[$(1)$] If $(X_{i+1},L_{i+1,p})$ satisfies Property $N_{p+1}$, then so
does $(X_{i},L_{i,p})$.
\item[$(2)$] $(X_2,L_{2,p})$ fails to satisfy Property
$N_{p+1}$.\\
\end{enumerate}
For $(1)$, consider the short exact sequence
\begin{equation*}
0 \rightarrow \mathcal{O}_{X_{i+1}} (-X_i) \rightarrow
\mathcal{O}_{X_{i+1}} \rightarrow \mathcal{O}_{X_{i}} \rightarrow
0
\end{equation*}
and thus
\begin{eqnarray*}
0 \rightarrow \wedge^m M_{i+1,p} \otimes L_{i+1,p} ^j \otimes
\mathcal{O}_{X_{i+1}} (-X_i) \rightarrow \wedge^m M_{i+1,p}
\otimes L_{i+1,p} ^j \\
\rightarrow \wedge^m M_{i+1,p} \otimes
L_{i+1,p} ^j \otimes \mathcal{O}_{X_{i}} \rightarrow 0
\end{eqnarray*}
In the derived cohomology long exact sequence
\begin{eqnarray*}
\cdots \rightarrow H^1 (X_{i+1},\wedge^m M_{i+1,p} \otimes
L_{i+1,p} ^j) \rightarrow H^1 (X_{i},\wedge^m M_{i+1,p} \otimes
L_{i+1,p} ^j \otimes \mathcal{O}_{X_{i}})\\
\rightarrow H^2 (X_{i+1},\wedge^m M_{i+1,p} \otimes L_{i+1,p} ^j
\otimes \mathcal{O}_{X_{i+1}} (-X_i)) \rightarrow \cdots,
\end{eqnarray*}
$H^1 (X_{i+1},\wedge^m M_{i+1,p} \otimes L_{i+1,p} ^j)=0$ for $1
\leq m \leq p+2$ and all $j \geq 1$ since $(X_{i+1},L_{i+1,p})$
satisfies Property $N_{p+1}$. Also $H^2 (X_{i+1},\wedge^m
M_{i+1,p} \otimes L_{i+1,p} ^j \otimes \mathcal{O}_{X_{i+1}}
(-X_i))=0$ for $1 \leq m \leq p+2$ and all $j \geq 1$. Therefore
$(X_i,L_{i,p})$ satisfies Property $N_{p+1}$. For $(2)$, let $C_0$
be the minimal section of $X_2$ and consider the short exact
sequence
\begin{equation*}
0 \rightarrow \mathcal{O}_{X_2} (-C_0) \rightarrow
\mathcal{O}_{X_2} \rightarrow \mathcal{O}_{C} \rightarrow 0.
\end{equation*}
Then by the same way one can prove that if $(X_2,L_{2,p})$
satisfies Property $N_{p+1}$, then so does $(C,L_C)$ where $L_C =
L_{2,p} |_{C_0}$ is the restriction of $L_{2,p}$ to $C_0$. But
$\mbox{deg}(L_C) = 2g+1+p$ and hence $L_C$ does not satisfy
Property $N_{p+1}$ since $C$ is a hyperelliptic curve. See
\cite{GL2}. Therefore $(X_2,L_{2,p})$ doesn't satisfy Property
$N_{p+1}$.
\end{proof}

\noindent {\bf Remark 5.1.} When $X$ is an anticanonical rational
surface, i.e., $X$ is a rational surface and $-K_X$ is effective, Mukai's conjecture is true.
See Theorem 1.23 and Example 1.15 in \cite{GP}. Indeed the
examples in which Mukai's conjecture is optimal is also classified
by their results. Let $A_i \in \mbox{Pic}X$ be ample line bundles
and assume that $K_X + A_1 + \cdots A_{p+4}$ fails to satisfy
Property $N_{p+1}$. Then $K_X ^2 =1$, $-K_X$ is ample and $A_1 =
\cdots = A_{p+4} = -K_X$. This can be easily checked
by using Theorem 1.3 and Proposition 1.10 in \cite{GP}. \qed \\


\begin{thebibliography}{00}
\bibitem{Butler} David C. Butler,
{\em Normal generation of vector bundles over a curve}, J. Differ.
Geom., 39 (1994), 1-34.

\bibitem{Eisenbud} David Eisenbud,
{\em The Geometry of Syzygies}, New York, forthcoming book(2004)

\bibitem{EGHP} D. Eisenbud, M. Green, K. Hulek and S. Popescu,
{\em Restriction linear syzygies: algebra and geometry},
math.AG/0404516.

\bibitem{EL} L. Ein and R. Lazarsfeld,
{\em Syzygies and Koszul cohomology of smooth projective varieties
of arbitrary dimension}, Inv. Math. 111 (1993), 51-67.

\bibitem{Green} M. Green,
{\em Koszul cohomology and the geometry of projective varieties
I}, J. Differ. Geom., 19 (1984), 125-171.

\bibitem{Green2} M. Green,
{\em Koszul cohomology and the geometry of projective varieties
II}, J. Differ. Geom., 20 (1984), 279-289.

\bibitem{GL} M. Green and R. Lazarsfeld,
{\em On the projective normality of complete linear series on an
algebraic curve}, Inv. Math. 83 (1986), 73-90.

\bibitem{GL2} M. Green and R. Lazarsfeld,
{\em Some results on the syzygies of finite sets and algebraic
curves}, Compositio Mathematica, 67 (1988), 301-314.

\bibitem{GP} F. J. Gallego and  B. P. Purnaprajna,
{\em Some results on rational surfaces and Fano varieties}, J.
Reine Angew. Math. 538(2001), 25-55

\bibitem{H} R. Hartshorne,
{\em Algebraic Geometry}, no. 52, Springer-Velag New York (1977)

\bibitem{KP} Sijong Kwak and Euisung Park,
{\em Some effects of property $N_p$ on the higher normality and
defining equations of nonlinearly normal varieties},
math.AG/0401026, to apprear in J. Reine Angew. Math.

\bibitem{OP} G. Ottaviani and R. Paoletti, {\em Syzygies of Veronese
embeddings}, Compositio Mathematica, 125 (2001), 31-37.

\bibitem{Park} Euisung Park, {\em On higher syzygies of ruled surfaces},
math.AG/0401100

\bibitem{Ru} Elena Rubei, {\em On syzygies of Segre embeddings},
Proc. Amer. Math. Soc, 130, no 12, (2002), 3483-3493.

\bibitem{Ru2} Elena Rubei, {\em A result on resolutions of
Veronese embeddings}, math.AG/0309102.
\end{thebibliography}
\end{document}